\documentclass[12pt,twoside]{article}
\usepackage[hmargin=0.8in,vmargin=0.9in]{geometry}
\usepackage{fancyhdr}
\usepackage{graphicx}
\usepackage{subfigure}
\usepackage{amssymb}
\usepackage{amsmath}
\usepackage{amsfonts}
\usepackage{theorem}
\usepackage{mathrsfs}
\usepackage{mathtools}
\usepackage{bm}
\usepackage{color}
\usepackage{setspace}
\usepackage{exscale}
\usepackage{relsize}
\usepackage{epstopdf}
\usepackage{float}
\usepackage{cite}

\usepackage{cleveref}

\usepackage{nicefrac}
\usepackage{setspace}

\usepackage{booktabs,multirow} 
\usepackage{array} 
\usepackage{paralist} 
\usepackage{verbatim} 
\usepackage{subfigure} 

\theoremstyle{plain}			

{\theorembodyfont{\rmfamily}\newtheorem{remark}{Remark}[section]}

\newenvironment{DA}{{\flushleft \bf Declarations:}}{}

\usepackage{tikz}
\usetikzlibrary{positioning}
\usepackage{xcolor}

\allowdisplaybreaks[1]

\numberwithin{equation}{section}
\numberwithin{figure}{section}
\numberwithin{table}{section}

\newcommand\eref[1]{(\ref{#1})}

\newcommand*\xbar[1]{%
  \hbox{%
    \vbox{%
      \hrule height 0.5pt 
      \kern0.4ex
      \hbox{%
        \kern-0.05em
        \ensuremath{#1}%
        \kern-0.00em
      }%
    }%
  }%
}

\newcommand{\bmF}{\bm{\mathcal{F}}}
\newcommand{\bmG}{\bm{\mathcal{G}}}

\newcommand{\bmK}{\bm{\mathcal{K}}}
\newcommand{\bmL}{\bm{\mathcal{L}}}
\newcommand{\mF}{\bm{F}}
\newcommand{\mK}{\bm{K}}
\newcommand{\mL}{\bm{L}}

\newcommand{\mG}{\bm{G}}

\newcommand{\mU}{\bm{U}}
\newcommand{\mE}{\bm{E}}

\newcommand{\mo}{\bm{0}}

\newcommand{\dt}{\Delta t}
\newcommand{\dx}{\Delta x}
\newcommand{\dy}{\Delta y}

\newcommand{\hf}{{\frac{1}{2}}}

\newcommand{\jph}{{j+\frac{1}{2}}}
\newcommand{\jmh}{{j-\frac{1}{2}}}
\newcommand{\kph}{{k+\frac{1}{2}}}
\newcommand{\kmh}{{k-\frac{1}{2}}}

\def\softd{{\leavevmode\setbox1=\hbox{d}%
          \hbox to 1.05\wd1{d\kern-0.4ex{\char039}\hss}}}

\title{New More Efficient A-WENO Schemes}
\author{Shaoshuai Chu\thanks{ Department of Mathematics, RWTH Aachen University, 52056 Aachen, Germany; Department of Mathematics and
Shenzhen International Center for Mathematics, Southern University of Science and Technology, Shenzhen, 518055, China;
{\tt chuss2019@mail.sustech.edu.cn}}, Alexander Kurganov\thanks{Department of Mathematics and Shenzhen International Center for Mathematics,
Southern University of Science and Technology, Shenzhen, 518055, China; {\tt alexander@sustech.edu.cn}}, and Ruixiao Xin\thanks{Department
of Mathematics, Southern University of Science and Technology, Shenzhen, 518055, China; {\tt 12331009@mail.sustech.edu.cn}}}

\begin{document}

\date{}
\maketitle
\begin{abstract}
We develop new more efficient A-WENO schemes for both hyperbolic systems of conservation laws and nonconservative hyperbolic systems.
The new schemes are a very simple modification of the existing A-WENO schemes: They are obtained by a more efficient evaluation
of the high-order correction terms. We conduct several numerical experiments to demonstrate the performance of the introduced schemes.
\end{abstract}



\section{Introduction}
This paper is focused on the development of new more efficient fifth-order finite-difference (FD) alternative weighted essentially
non-oscillatory (A-WENO) scheme for both conservative and nonconservative hyperbolic systems of nonlinear PDEs.

In the one- and two-dimensional cases, hyperbolic systems of conservation laws read as
\begin{equation}
\bm U_t+\bm F(\bm U)_x=\mo,
\label{1.1}
\end{equation}
\begin{equation}
\bm U_t+\bm F(\bm U)_x+\bm G(\bm U)_y=\mo,
\label{1.1a}
\end{equation}
respectively. In \eref{1.1}--\eref{1.1a}, $x$ and $y$ are the spatial variables, $t$ is the time, $\mU\in\mathbb R^d$ is a vector of
unknowns, $\mF:\mathbb R^d\to\mathbb R^d$ and $\mG:\mathbb R^d\to\mathbb R^d$ are nonlinear flux functions. The A-WENO schemes for
\eref{1.1} and \eref{1.1a}, introduced in \cite{Jiang13} (also see
\cite{CCK23_Adaptive,Don20,liu17,liu16,WDGK_20,wang18,WDKL,CKX22}), employ finite-volume (FV) numerical fluxes without any flux splitting,
and use high-order correction terms to ensure a high-order spatial accuracy of the resulting scheme. For instance, the original A-WENO
scheme from \cite{Jiang13} employs an upwind flux, the schemes in \cite{Don20,liu17,liu16,WDGK_20,wang18} rely on the simplest---yet very
robust---central (local Lax-Friedrichs, Rusanov) flux or its adaptive version, and the schemes in \cite{WDKL,CKX22,CCK23_Adaptive} are based
on central-upwind fluxes. Examples of many other FV numerical fluxes can be found in, e.g., \cite{Hesthaven18,Leveque02,KLR20,Tor}.

In this paper, we introduce a novel approach for computing high-order correction terms in A-WENO schemes. Traditional methods require
evaluating the flux functions $\mF(\mU)$ and $\mG(\mU)$ at each grid point during every time step, which can be computationally expensive.
The proposed method leverages precomputed numerical fluxes, allowing for time savings while maintaining both accuracy and stability without
risking any additional spurious oscillations.

The proposed modification can also be applied to the nonconservative hyperbolic systems, which in the one- and two-dimensional
cases read as
\begin{equation}
\bm U_t+\bm F(\bm U)_x=C(\bm U)\bm U_x,
\label{1.2}
\end{equation}
and 
\begin{equation}
\bm U_t+\bm F(\bm U)_x+\bm G(\bm U)_y=B(\bm U)\bm U_x+C(\bm U)\bm U_y,
\label{1.2b}
\end{equation}
respectively, with $B(\mU)\in\mathbb R^{d\times d}$ and $C(\mU)\in\mathbb R^{d\times d}$. The presence of the nonconservative products on
the right-hand side (RHS) of \eref{1.2} and \eref{1.2b} introduces more challenges for studying \eref{1.2} and
\eref{1.2b} both theoretically and numerically. Weak solutions of \eref{1.2} and \eref{1.2b} can be understood in the sense of
Borel measures as it was done in \cite{maso1995,lefloch02,lefloch2012}, but not in the sense of distributions when $\mU$ is discontinuous.
The concept of Borel measure solutions was utilized to develop path-conservative schemes. A wide of variety of path-conservative methods
have been introduced; see, e.g. \cite{pares09,castro17,CP20,PGCCMP,SGBNP,BDGI,Cha,CKM,KLX_21,CKLX_22,CKLZ_23} and references therein.
Recently, the second-order path-conservative central-upwind scheme from \cite{CKM} has been extended to the fifth order of
accuracy in the framework of A-WENO schemes; see, e.g., \cite{Chu21,CKMZ23}. In these schemes, we first rewrite the systems \eref{1.2} and
\eref{1.2b} in the following quasi-conservative forms:
\begin{equation}
\bm U_t+\bm K(\bm U)_x=\bm0,
\label{1.2a}
\end{equation}
where
\begin{equation}
\bm K(\bm U)=\bm F(\bm U)-\bm R(\bm U),\quad\bm R(\bm U):=\int\limits^x_{\widehat x}B(\bm U){\bm U}_\xi(\xi,t)\,{\rm d}\xi,
\label{1.2b1}
\end{equation}
and 
\begin{equation}
\bm U_t+\bm K(\bm U)_x+\bm L(\bm U)_y=\bm0,
\label{1.2c}
\end{equation}
where
\begin{equation}
\begin{aligned}
&\bm K(\bm U)=\bm F(\bm U)-\bm R^x(\bm U),\quad\bm R^x(\bm U):=\int\limits^x_{\widehat x}B(\bm U){\bm U}_\xi(\xi,y,t)\,{\rm d}\xi,\\
&\bm L(\bm U)=\bm G(\bm U)-\bm R^y(\bm U),\quad\bm R^y(\bm U):=\int\limits^y_{\widehat y}C(\bm U){\bm U}_\eta(x,\eta,t)\,{\rm d}\eta,
\end{aligned}
\label{1.2d}
\end{equation}
respectively. In \eref{1.2b1} and \eref{1.2d}, $\widehat x$ and $\widehat y$ are arbitrary numbers. High-order correction terms are then
computed using the point values of the global fluxes $\bm K$ and $\bm L$. As in the conservative case, the point values of
$\bm K$ and $\bm L$ need to be computed at every grid point at every time step. In this paper, we propose a more efficient
alternative: to use the already computed numerical flux values in the computation of the high-order correction terms, which leads to new
more efficient A-WENO schemes for nonconservative hyperbolic systems. The proposed A-WENO schemes are realized in the framework of flux
globalization based well-balanced (WB) path-conservative schemes, which were introduced in \cite{KLX_21} and later applied to a variety of
shallow water models in \cite{CKLX_22,CKLZ_23,CKN23}.

The rest of the paper is organized as follows. In \S\ref{sec2.1}, we briefly describe the existing one-dimensional (1-D)  A-WENO
schemes and then in \S\ref{sec2.2}, we introduce new 1-D A-WENO schemes using an alternative, more efficient way of computing
the high-order correction terms. In \S\ref{sec3.1}, we describe 1-D flux globalization based WB path-conservative A-WENO
schemes, and then introduce new more efficient 1-D flux globalization based WB path-conservative A-WENO schemes in
\S\ref{sec3.2}. In \S\ref{sec4} and \S\ref{sec5}, we extend the 1-D A-WENO schemes introduced in \S\ref{sec2} and \S\ref{sec3}
to the two-dimensional (2-D) cases.  Finally, in \S\ref{sec6}, we test the proposed schemes on a number of conservative and nonconservative
examples. The obtained numerical results demonstrate that the new A-WENO schemes are clearly more efficient than the old ones, while
neither the accuracy nor the quality of resolution is affected by switching to the new approach.

\section{Fifth-Order A-WENO Schemes for \eref{1.1}}\label{sec2}
\subsection{Existing A-WENO Schemes}\label{sec2.1}
We assume that the computational domain is covered with uniform cells $C_j:=[x_\jmh,x_\jph]$ of size $x_\jph-x_\jmh\equiv\dx$ centered at
$x_j=(x_\jmh+x_\jph)/2$, $j=1,\ldots,N_x$. We also assume that at a certain time $t\ge0$, the point values of the computed solution,
$\mU_j(t)$, are available (in the rest of the paper, we will suppress the time-dependence of all of the indexed quantities for the sake of
brevity). In the framework of semi-discrete A-WENO schemes, $\mU_j$ are evolved in time by solving the following system of ODEs:
\begin{equation}
\frac{{\rm d}\mU_j}{{\rm d}t}=-\frac{\bmF_\jph-\bmF_\jmh}{\dx},
\label{2.1a}
\end{equation}
where $\bmF_\jph$ are the fifth-order numerical fluxes, computed by
\begin{equation}
\bmF_\jph=\bm{{\cal F}}^{\rm FV}_\jph-\frac{(\dx)^2}{24}(\bm F_{xx})_\jph+\frac{7(\dx)^4}{5760}(\bm F_{xxxx})_\jph.
\label{2.1b}
\end{equation}
Here, $\bmF^{\rm FV}_\jph=\bmF^{\rm FV}_\jph\big(\mU^-_\jph,\mU^+_\jph\big)$ stand for FV numerical fluxes, which depends on $\mU^-_\jph$
and $\mU^+_\jph$, which are the left- and right-sided values of $\mU$ at $x=x_\jph$ computed using the WENO-Z interpolation procedure
introduced in \cite{Borges08,Castro11,Don13}. In this paper, we have utilized one of the simplest FV numerical fluxes---the central (local
Lax-Friedrichs, Rusanov) fluxes (see \cite{KTcl,Rus61}):
\begin{equation*}
\bm{{\cal F}}^{\rm FV}_\jph=\frac{\bm F(\mU^-_\jph)+\bm F(\mU^+_\jph) }{2}-\frac{a_\jph}{2}\big( \bm U^+_\jph-\bm U^-_\jph \big),
\end{equation*}
where the local speed of propagation $a_\jph$ is estimated using the eigenvalues $\lambda_1(A)\le\ldots\le\lambda_d(A)$ of the Jacobian
$A=\frac{\partial\mF}{\partial\mU}$:
\begin{equation*}
a_\jph=\max_j\Big\{\big|\lambda_1\big(A(\mU^-_\jph)\big)\big|,\,\big|\lambda_1\big(A(\mU^+_\jph)\big)\big|,\,
\big|\lambda_d\big(A(\mU^-_\jph)\big)\big|,\,\big|\lambda_d\big(A(\mU^+_\jph)\big)\big|\Big\}.
\end{equation*}

Finally, the higher-order correction terms $(\bm F_{xx})_\jph$ and $({\bm F_{xxxx}})_\jph$ are the second- and fourth-order spatial
numerical derivatives of $\bm F$ at the cell interface $x=x_\jph$, which are computed using the following FD approximations:
\begin{equation}
\begin{aligned}
&(\mF_{xx})_\jph=\frac{1}{48(\dx)^2}\left(-5\mF_{j-2}+39\mF_{j-1}-34\mF_j-34\mF_{j+1}+39\mF_{j+2}-5\mF_{j+3}\right),\\
&(\mF_{xxxx})_\jph=\frac{1}{2(\dx)^4}\left(\mF_{j-2}-3\mF_{j-1}+2\mF_j+2\mF_{j+1}-3\mF_{j+2}+\mF_{j+3}\right),
\end{aligned}
\label{2.1}
\end{equation}
where $\mF_j:=\mF(\mU_j)$.

\subsection{New A-WENO Schemes}\label{sec2.2}
In this section, we provide an alternative, more efficient way of computing the high-order correction terms, which are now
obtained by applying fourth- and second-order central difference formulae to evaluate $(\bm F_{xx})_\jph$ and $({\bm F_{xxxx}})_\jph$,
respectively:
\begin{equation}
\begin{aligned}
&(\mF_{xx})_\jph=\frac{1}{12(\dx)^2}\Big[-\bm{{\cal F}}^{\rm FV}_{j-\frac{3}{2}}+16\bm{{\cal F}}^{\rm FV}_\jmh-30\bm{{\cal F}}^{\rm FV}_\jph
+16\bm{{\cal F}}^{\rm FV}_{j+\frac{3}{2}}-\bm{{\cal F}}^{\rm FV}_{j+\frac{5}{2}}\Big],\\
&(\mF_{xxxx})_\jph=\frac{1}{(\dx)^4}\Big[\bm{{\cal F}}^{\rm FV}_{j-\frac{3}{2}}-4\bm{{\cal F}}^{\rm FV}_\jmh+6\bm{{\cal F}}^{\rm FV}_\jph-
4\bm{{\cal F}}^{\rm FV}_{j+\frac{3}{2}}+\bm{{\cal F}}^{\rm FV}_{j+\frac{5}{2}}\Big].
\end{aligned}
\label{2.3}
\end{equation}
Notice that the main difference between \eref{2.3} and \eref{2.1} is that in \eref{2.3}, we use FD approximations of $(\bm F_{xx})_\jph$ and
$({\bm F_{xxxx}})_\jph$ terms with the help of the FV numerical fluxes $\bm{{\cal F}}^{\rm FV}_\jph$ instead of the point values $\bm F_j$.
Clearly, the new A-WENO schemes are less computationally expensive than the original ones as they do not require computing the point values
$\mF_j$, while the values $\bm{{\cal F}}^{\rm FV}_\jph$ have to be computed in \eref{2.1b} and then can be stored. At the same time, the new
A-WENO schemes are still fifth-order accurate and they can achieve the same high resolution as the original A-WENO schemes as demonstrated
in \S\ref{sec6}, where we present several numerical examples.
\begin{remark}
The new high-order correction terms can be extended to higher than the fifth order of accuracy via higher-order Taylor
expansions; see, e.g., \cite{Gao20} for details on the seventh- and ninth-order A-WENO schemes.
\end{remark}

\begin{remark}
It should be pointed out that the use of higher-order correction terms \eref{2.1} may lead to the appearance of spurious oscillations; see, e.g., \cite{BBSK2024a,BBSK2024b,BBSK2024c}. We would like to stress that the use of the proposed higher-order correction terms \eref{2.3} does not help to suppress these oscillations. The way to reduce them was proposed in \cite{BBSK2024a,BBSK2024b,BBSK2024c}.
\end{remark}

\section{Flux Globalization Based Well-Balanced Path-Conservative A-WENO Schemes for \eref{1.2}}\label{sec3}
In this section, we extend the A-WENO scheme presented in \S\ref{sec2} to the nonconservative systems \eref{1.2} via a flux globalization
based framework.

\subsection{Flux Globalization Based A-WENO Schemes}\label{sec3.1}
We begin with the extension of A-WENO schemes to the nonconservative systems \eref{1.2}, written in the quasi-conservative form \eref{1.2a}--\eref{1.2b1}.
A direct generalization of the semi-discrete A-WENO scheme \eref{2.1a}--\eref{2.1b} reads as
\begin{equation}
\frac{{\rm d}\bm U_j}{{\rm d}t}=-\frac{\bmK_\jph-\bmK_\jmh}{\dx},
\label{3.1}
\end{equation}
where $\bmK_\jph$ are the fifth-order numerical fluxes:
\begin{equation*}
\bmK_\jph=\bm{{\cal K}}^{\rm FV}_\jph-\frac{(\dx)^2}{24}(\bm K_{xx})_\jph+\frac{7(\dx)^4}{5760}(\bm K_{xxxx})_\jph.
\end{equation*}
Here, $\bm{{\cal K}}^{\rm FV}_\jph=\bm{{\cal K}}^{\rm FV}_\jph\big(\mU^-_\jph,\mU^+_\jph\big)$ stand for FV numerical fluxes and in this
paper, we use the central (local Lax-Friedrichs, Rusanov) fluxes
\begin{equation}
\bm{{\cal K}}^{\rm FV}_\jph=\frac{\bm K\big(\mU^-_\jph\big)+\bm K\big(\mU^+_\jph\big)}{2}-\frac{a_\jph}{2}
\big(\bm U^+_\jph-\bm U^-_\jph\big)
\label{3.2}
\end{equation}
with the local speeds of propagation estimated using the eigenvalues $\lambda_1({\cal A})\le\ldots\le\lambda_d({\cal A})$ of the matrix
${\cal{A}}=\frac{\partial\bm F}{\partial\bm U}(\bm U)-B(\bm U)$:
\begin{equation*}
a_\jph=\max_j\Big\{\big|\lambda_1\big({\cal A}(\mU^-_\jph)\big)\big|,\,\big|\lambda_1\big({\cal A}(\mU^+_\jph)\big)\big|,\,
\big|\lambda_d\big({\cal A}(\mU^-_\jph)\big)\big|,\,\big|\lambda_d\big({\cal A}(\mU^+_\jph)\big)\big|\Big\}.
\end{equation*}

The global fluxes $\bm K^\pm_\jph=\bm K(\mU_\jph^\pm)$ in \eref{3.2} are to be computed in such a way to ensure both WB and
path-conservative properties of the resulting scheme. To this end, we follow the approach introduced in \cite{KLX_21}. We first note that
\eref{1.2} admits steady-state solutions satisfying
\begin{equation}
\bm F(\bm U)_x-B(\bm U)\bm U_x=\bm0.
\label{3.1a}
\end{equation}
As it was shown in \cite{KLX_21} (also see \cite{CKLX_22,CKLZ_23,CKN23}), for many particular nonconservative systems the relation
\eref{3.1a} can be written as
\begin{equation}
\bm F(\bm U)_x-B(\bm U)\bm U_x=M(\bm U)\bm E(\bm U)_x=\bm0,
\label{3.1b}
\end{equation}
where $M\in\mathbb R^{d\times d}$ and $\bm E$ is the vector of equilibrium variables. If $M$ is invertible (which is the case in many
practically interesting cases), the steady states satisfy the relation
\begin{equation*}
\bm E(\bm U)\equiv\mbox{\bf Const}.
\end{equation*}
Taking this into account, we evaluate $\mK^\pm_\jph$ as follows:
\begin{equation}
\bm K^\pm_\jph=\bm F^\pm_\jph-\bm R^\pm_\jph,
\label{3.3a}
\end{equation}
where $\bm F^\pm_\jph:=\bm F(\bm U^\pm_\jph)$, and the point values $\bm U^\pm_\jph$ are to be obtained using a reconstruction of the
equilibrium variables $\bm E$ to ensure the WB property. To this end, we first compute the point values $\mE_j:=\mE(\mU_j)$ and then apply
the WENO-Z interpolation to evaluate $\bm E^\pm_\jph$. The corresponding values $\bm U^\pm_\jph$ are then computed by solving the (nonlinear
systems of) equations
\begin{equation*}
\bm E(\bm U^+_\jph)=\bm E^+_\jph\quad\mbox{and}\quad\bm E(\bm U^-_\jph)=\bm E^-_\jph
\end{equation*}
for $\bm U^+_\jph$ and $\bm U^-_\jph$, respectively. The point values of the global variable $\bm R^\pm_\jph$ in \eref{3.3a} are computed in
a path-conservative way.  We first select $\widehat x=x_\hf$, set $\bm R^-_\hf:=\bm0$, evaluate
\begin{equation}
\bm R^+_\hf=\bm B_{\bm\Psi,\hf},
\label{3.5}
\end{equation}
and then recursively obtain
\begin{equation}
\bm R^-_\jph=\bm R^+_\jmh+\bm B_j,\quad\bm R^+_\jph=\bm R^-_\jph+\bm B_{\bm\Psi,\jph},\quad j=1,\ldots,N_x.
\label{3.6}
\end{equation}
In \eref{3.5} and \eref{3.6}, $\bm B_j$ and $\bm B_{\bm\Psi,\jph}$ are obtained using \eref{3.1b} and proper quadratures:
\begin{equation}
\begin{aligned}
&\bm B_j=\int\limits_{C_j}B\big(\bm U\big)\bm U_x\,{\rm d}x=\bm F_\jph^--\bm F_\jmh^+-\int\limits_{C_j}M\big(\bm U\big)\bm E(\bm U)_x\,{\rm d}x,\\
&\bm B_{\bm\Psi,\jph}=\int\limits^1_0B\big(\bm\Psi_\jph(s)\big)\bm\Psi'_\jph(s)\,{\rm d}s\\
&\hspace{1.3cm}=\bm F_\jph^+-\bm F_\jph^--\hf\left[M\big(\bm U_\jph^+\big)+M\big(\bm U_\jph^-\big)\right]
\big(\bm E_\jph^+-\bm E_\jph^-\big),
\label{3.8}
\end{aligned}
\end{equation}
where $\bm\Psi_\jph(s):=\bm\Psi_\jph(s;\bm U^-_\jph,\bm U^+_\jph)$ is a sufficiently smooth path connecting the left- and right-sided states
$\bm U^-_\jph=\bm U\big(\bm E^-_\jph\big)$ and $\bm U^+_\jph=\bm U\big(\bm E^+_\jph\big)$. We refer the readers to \cite{CKX23_WB} for more
details on the derivation of \eref{3.8}.

In order to ensure the designed flux globalization based  WB path-conservative A-WENO scheme is fifth-order accurate, the integral $\bm B_j$
in \eref{3.8} needs to be evaluated using at least a fifth-order quadrature. In this paper, we use the Newton-Cotes quadrature introduced
in \cite[(4.4)]{Chu21}.

Finally, $(\bm K_{xx})_\jph$ and $({\bm K_{xxxx}})_\jph$ are approximations of the second- and fourth-order spatial derivatives of $\bm K$
at $x=x_\jph$, which can be computed using the same FD approximations as in \eref{2.1}, namely, by
\begin{equation*}
\begin{aligned}
&(\mK_{xx})_\jph=\frac{1}{48(\dx)^2}\left(-5\mK_{j-2}+39\mK_{j-1}-34\mK_j-34\mK_{j+1}+39\mK_{j+2}-5\mK_{j+3}\right),\\
&(\mK_{xxxx})_\jph=\frac{1}{2(\dx)^4}\left(\mK_{j-2}-3\mK_{j-1}+2\mK_j+2\mK_{j+1}-3\mK_{j+2}+\mK_{j+3}\right),
\end{aligned}
\end{equation*}
where $\mK_j=\mF(\mU_j)-\bm R_j(\bm U)$ and
\begin{equation}
{\bm R_j}(\bm U)=\int\limits_{\widehat x}^{x_j}B({\bm U(\xi,t)})\bm U_\xi(\xi,t)\,{\rm d}\xi.
\label{3.16}
\end{equation}
In order to obtain a fifth-order numerical scheme, the integral in \eref{3.16} needs to be evaluated within the fifth order of accuracy.
We refer the reader to, e.g., \cite{CKMZ23} for the fifth-order Newton-Cotes quadrature, which should be implemented in a
recursive way.

\subsection{More Efficient Flux Globalization Based A-WENO Schemes}\label{sec3.2}
In this section, we extend a new, more efficient A-WENO scheme from \S\ref{sec2.2} to the nonconservative system \eref{1.2}. We now compute
the high-order correction terms $(\bm K_{xx})_\jph$ and $({\bm K_{xxxx}})_\jph$ with the help of the same central difference formulae (as was used in \S\ref{sec2.2}) to the numerical fluxes:
\begin{equation}
\begin{aligned}
&(\mK_{xx})_\jph=\frac{1}{12(\dx)^2}\Big[-\bm{{\cal K}}^{\rm FV}_{j-\frac{3}{2}}+16\bm{{\cal K}}^{\rm FV}_\jmh-30\bm{{\cal K}}^{\rm FV}_\jph
+16\bm{{\cal K}}^{\rm FV}_{j+\frac{3}{2}}-\bm{{\cal K}}^{\rm FV}_{j+\frac{5}{2}}\Big],\\
&(\mK_{xxxx})_\jph=\frac{1}{(\dx)^4}\Big[\bm{{\cal K}}^{\rm FV}_{j-\frac{3}{2}}-4\bm{{\cal K}}^{\rm FV}_\jmh+6\bm{{\cal K}}^{\rm FV}_\jph-
4\bm{{\cal K}}^{\rm FV}_{j+\frac{3}{2}}+\bm{{\cal K}}^{\rm FV}_{j+\frac{5}{2}}\Big].
\end{aligned}
\label{3.10f}
\end{equation}

Notice that here we do not need to evaluate the point values $\mK_j$ and, in particular, we do not need to evaluate the integrals in
\eref{3.16}. Therefore, the resulting new version of the flux globalization based WB path-conservative A-WENO scheme is substantially less
computationally expensive than the scheme from \S\ref{sec3.1}.
\begin{remark}
As in the conservative case, the presented flux globalization based WB path-conservative A-WENO scheme can be straightforwardly extended to
higher than the fifth order of accuracy via higher-order Taylor expansions.
\end{remark}

\section{Fifth-Order A-WENO Schemes for \eref{1.1a}}\label{sec4}
\subsection{Existing Two-Dimensional A-WENO Schemes}\label{sec4.1}
We assume that the computational domain is covered with uniform cells $C_{j,k}:=[x_\jmh,x_\jph]\times[y_\kmh,y_\kph]$ with
$x_\jph-x_\jmh\equiv\dx$ and $y_\kph-y_\kmh\equiv\dy$ centered at $(x_j,y_k)$, $j=1,\ldots,N_x,k=1,\ldots,N_y$, with $x_j=(x_\jmh+x_\jph)/2$
and $y_k=(y_\kmh+y_\kph)/2$. We also assume that at a certain time $t\ge0$, the point values of the computed solution, $\mU_{j,k}(t)$, are
available. In the framework of semi-discrete A-WENO schemes, $\mU_{j,k}$ are evolved in time by solving the following system of ODEs:
\begin{equation}
\frac{{\rm d}\mU_{j,k}}{{\rm d}t}=-\frac{\bmF_{\jph,k}-\bmF_{\jmh,k}}{\dx}-\frac{\bmG_{j,\kph}-\bmG_{j,\kmh}}{\dy},
\label{4.1a}
\end{equation}
where $\bmF_{\jph,k}$ and $\bmG_{j,\kph}$ are the fifth-order numerical fluxes, computed by
\begin{equation}
\begin{aligned}
& \bmF_{\jph,k}=\bm{{\cal F}}^{\rm FV}_{\jph,k}-\frac{(\dx)^2}{24}(\bm F_{xx})_{\jph,k}+\frac{7(\dx)^4}{5760}(\bm F_{xxxx})_{\jph,k},\\
& \bmK_{j,\kph}=\bm{{\cal G}}^{\rm FV}_{j,\kph}-\frac{(\dy)^2}{24}(\bm G_{yy})_{j,\kph}+\frac{7(\dy)^4}{5760}(\bm G_{yyyy})_{j,\kph}.
\end{aligned}
\label{4.1b}
\end{equation}
Here, $\bmF^{\rm FV}_{\jph,k}=\bmF^{\rm FV}_{\jph,k}\big(\mU^-_{\jph,k},\mU^+_{\jph,k}\big)$ and
$\bmG^{\rm FV}_{j,\kph}=\bmG^{\rm FV}_{j,\kph}\big(\mU^-_{j,\kph},\mU^+_{j,\kph}\big)$ stand for central (local Lax-Friedrichs, Rusanov)
fluxes (see \cite{KTcl,Rus61}):
\begin{equation*}
\begin{aligned}
&\bm{{\cal F}}^{\rm FV}_{\jph,k}=\frac{\bm F(\mU^-_{\jph,k})+\bm F(\mU^+_{\jph,k}) }{2}-
\frac{a_{\jph,k}}{2}\big( \bm U^+_{\jph,k}-\bm U^-_{\jph,k} \big),\\
&\bm{{\cal G}}^{\rm FV}_{j,\kph}=\frac{\bm G(\mU^-_{j,\kph})+\bm G(\mU^+_{j,\kph}) }{2}-
\frac{a_{j,\kph}}{2}\big( \bm U^+_{j,\kph}-\bm U^-_{j,\kph} \big),
\end{aligned}
\end{equation*}
where $\mU^\pm_{\jph,k}$ and $\mU^\pm_{j,\kph}$ are the one-sided point values of $\mU$ at the cell interfaces $(x_\jph,y_k)$ and
$(x_j,y_\kph)$, respectively, which are computed using the 1-D WENO-Z interpolations performed in the $x$- and $y$-directions, respectively.
The local speed of propagation in the $x$- and $y$-directions, $a_{\jph,k}$ and $a_{j,\kph}$, are estimated using the smallest and largest
eigenvalues of the Jacobians $A=\frac{\partial\mF}{\partial\mU}$ and $B=\frac{\partial\mG}{\partial\mU}$:
\begin{equation*}
\begin{aligned}
&a_{\jph,k}=\max_{j,k}\Big\{\big|\lambda_1\big(A(\mU^-_{\jph,k})\big)\big|,\,\big|\lambda_1\big(A(\mU^+_{\jph,k})\big)\big|,\,
\big|\lambda_d\big(A(\mU^-_{\jph,k})\big)\big|,\,\big|\lambda_d\big(A(\mU^+_{\jph,k})\big)\big|\Big\},\\
&a_{j,\kph}=\max_{j,k}\Big\{\big|\lambda_1\big(B(\mU^-_{j,\kph})\big)\big|,\,\big|\lambda_1\big(B(\mU^+_{j,\kph})\big)\big|,\,
\big|\lambda_d\big(B(\mU^-_{j,\kph})\big)\big|,\,\big|\lambda_d\big(B(\mU^+_{j,\kph})\big)\big|\Big\}.
\end{aligned}
\end{equation*}

Finally, the higher-order correction terms $(\bm F_{xx})_{\jph,k}$, $({\bm F_{xxxx}})_{\jph,k}$, $(\bm G_{yy})_{j,\kph}$, and
$({\bm G_{yyyy}})_{j,\kph}$ are the second- and fourth-order spatial numerical derivatives of $\bm F$ and $\bm G$ at the cell interfaces
$(x_\jph,y_k)$ and $(x_j,y_\kph)$, respectively. They are computed using the following FD approximations:
\begin{equation}
\begin{aligned}
&(\mF_{xx})_{\jph,k}=\frac{1}{48(\dx)^2}\left(-5\mF_{j-2,k}+39\mF_{j-1,k}-34\mF_{j,k}-34\mF_{j+1,k}+39\mF_{j+2,k}-5\mF_{j+3,k}\right),\\
&(\mF_{xxxx})_{\jph,k}=\frac{1}{2(\dx)^4}\left(\mF_{j-2,k}-3\mF_{j-1,k}+2\mF_{j,k}+2\mF_{j+1,k}-3\mF_{j+2,k}+\mF_{j+3,k}\right),\\
&(\mG_{yy})_{j,\kph}=\frac{1}{48(\dy)^2}\left(-5\mG_{j,k-2}+39\mG_{j,k-1}-34\mG_{j,k}-34\mG_{j,k+1}+39\mG_{j,k+2}-5\mG_{j,k+3}\right),\\
&(\mG_{yyyy})_{j,\kph}=\frac{1}{2(\dy)^4}\left(\mG_{j,k-2}-3\mG_{j,k-1}+2\mG_{j,k}+2\mG_{j,k+1}-3\mG_{j,k+2}+\mG_{j,k+3}\right),
\end{aligned}
\label{4.1}
\end{equation}
where $\mF_{j,k}:=\mF(\mU_{j,k})$ and $\mG_{j,k}:=\mG(\mU_{j,k})$.

\subsection{New Two-Dimensional A-WENO Schemes}\label{sec4.2}
In this section, we provide an alternative, more efficient way of computing the high-order correction terms$(\bm F_{xx})_{\jph,k}$,
$({\bm F_{xxxx}})_{\jph,k}$, $(\bm G_{yy})_{j,\kph}$, and $({\bm G_{yyyy}})_{j,\kph}$, which are obtained using the same formulae as in
\S\ref{sec2.2} and \S\ref{sec3.2}:
\begin{equation}
\begin{aligned}
&(\mF_{xx})_{\jph,k}=\frac{1}{12(\dx)^2}\Big[-\bm{{\cal F}}^{\rm FV}_{j-\frac{3}{2},k}+16\bm{{\cal F}}^{\rm FV}_{\jmh,k}-
30\bm{{\cal F}}^{\rm FV}_{\jph,k}+16\bm{{\cal F}}^{\rm FV}_{j+\frac{3}{2},k}-\bm{{\cal F}}^{\rm FV}_{j+\frac{5}{2},k}\Big],\\
&(\mF_{xxxx})_{\jph,k}=\frac{1}{(\dx)^4}\Big[\bm{{\cal F}}^{\rm FV}_{j-\frac{3}{2},k}-4\bm{{\cal F}}^{\rm FV}_{\jmh,k}+
6\bm{{\cal F}}^{\rm FV}_{\jph,k}-4\bm{{\cal F}}^{\rm FV}_{j+\frac{3}{2},k}+\bm{{\cal F}}^{\rm FV}_{j+\frac{5}{2},k}\Big],\\
&(\mG_{yy})_{j,\kph}=\frac{1}{12(\dy)^2}\Big[-\bm{{\cal G}}^{\rm FV}_{j,k-\frac{3}{2}}+16\bm{{\cal G}}^{\rm FV}_{j,\kmh}-
30\bm{{\cal G}}^{\rm FV}_{j,\kph}+16\bm{{\cal G}}^{\rm FV}_{j,k+\frac{3}{2}}-\bm{{\cal G}}^{\rm FV}_{j,k+\frac{5}{2}}\Big],\\
&(\mG_{yyyy})_{j,\kph}=\frac{1}{(\dy)^4}\Big[\bm{{\cal G}}^{\rm FV}_{j,k-\frac{3}{2}}-4\bm{{\cal G}}^{\rm FV}_{j,\kmh}+
6\bm{{\cal G}}^{\rm FV}_{j,\kph}-4\bm{{\cal G}}^{\rm FV}_{j,k+\frac{3}{2}}+\bm{{\cal G}}^{\rm FV}_{j,k+\frac{5}{2}}\Big].
\end{aligned}
\label{4.4f}
\end{equation}
\begin{remark}
As in the 1-D case, the new high-order correction terms can be extended to higher than the fifth order of accuracy via higher-order Taylor
expansions.
\end{remark}

\section{Flux Globalization Based Well-Balanced Path-Conservative A-WENO Schemes for \eref{1.2b}}\label{sec5}
In this section, we extend the 2-D A-WENO scheme presented in \S\ref{sec4} to the nonconservative systems \eref{1.2b} within the flux
globalization based framework.

\subsection{Two-Dimensional Flux Globalization Based A-WENO Schemes}\label{sec5.1}
We first rewrite the nonconservative systems \eref{1.2b} in the quasi-conservative form \eref{1.2c}--\eref{1.2d}. A direct generalization of
the semi-discrete A-WENO scheme \eref{4.1a}--\eref{4.1b} reads as
\begin{equation}
\frac{{\rm d}\bm U_{j,k}}{{\rm d}t}=-\frac{\bmK_{\jph,k}-\bmK_{\jmh,k}}{\dx}-\frac{\bmL_{j,\kph}-\bmK_{j,\kmh}}{\dy},
\label{5.1}
\end{equation}
where $\bmK_{\jph,k}$ and $\bmL_{j,\kph}$ are the fifth-order numerical fluxes in the $x$- and $y$-directions:
\begin{equation*}
\begin{aligned}
&\bmK_{\jph,k}=\bm{{\cal K}}^{\rm FV}_{\jph,k}-\frac{(\dx)^2}{24}(\bm K_{xx})_{\jph,k}+\frac{7(\dx)^4}{5760}(\bm K_{xxxx})_{\jph,k},\\
&\bmL_{j,\kph}=\bm{{\cal L}}^{\rm FV}_{j,\kph}-\frac{(\dy)^2}{24}(\bm L_{yy})_{j,\kph}+\frac{7(\dx)^4}{5760}(\bm L_{yyyy})_{j,\kph}.
\end{aligned}
\end{equation*}
Here, $\bm{{\cal K}}^{\rm FV}_{\jph,k}=\bm{{\cal K}}^{\rm FV}_{\jph,k}\big(\mU^-_{\jph,k},\mU^+_{\jph,k}\big)$ and
$\bm{{\cal L}}^{\rm FV}_{j,\kph}=\bm{{\cal L}}^{\rm FV}_{j,\kph}\big(\mU^-_{j,\kph},\mU^+_{j,\kph}\big)$ stand for FV numerical fluxes and
we use the central (local Lax-Friedrichs, Rusanov) fluxes
\begin{equation}
\begin{aligned}
&\bm{{\cal K}}^{\rm FV}_{\jph,k}=\frac{\bm K\big(\mU^-_{\jph,k}\big)+\bm K\big(\mU^+_{\jph,k}\big)}{2}-\frac{a_{\jph,k}}{2}
\big(\bm U^+_{\jph,k}-\bm U^-_{\jph,k}\big),\\
&\bm{{\cal L}}^{\rm FV}_{j,\kph}=\frac{\bm L\big(\mU^-_{j,\kph}\big)+\bm L\big(\mU^+_{j,\kph}\big)}{2}-\frac{a_{j,\kph}}{2}
\big(\bm U^+_{j,\kph}-\bm U^-_{j,\kph}\big),
\label{5.2}
\end{aligned}
\end{equation}
with the local speeds of propagation estimated using the smallest and largest eigenvalues of
${\cal A}=\frac{\partial\bm F}{\partial\bm U}(\bm U)-B(\bm U)$ and ${\cal B}=\frac{\partial\bm G}{\partial\bm U}(\bm U)-C(\bm U)$:
\begin{equation*}
\begin{aligned}
&a_{\jph,k}=\max_{j,k}\Big\{\big|\lambda_1\big({\cal A}(\mU^-_{\jph,k})\big)\big|,\,\big|\lambda_1\big({\cal A}(\mU^+_{\jph,k})\big)\big|,\,
\big|\lambda_d\big({\cal A}(\mU^-_{\jph,k})\big)\big|,\,\big|\lambda_d\big({\cal A}(\mU^+_{\jph,k})\big)\big|\Big\},\\
&a_{j,\kph}=\max_{j,k}\Big\{\big|\lambda_1\big({\cal B}(\mU^-_{j,\kph})\big)\big|,\,\big|\lambda_1\big({\cal B}(\mU^+_{j,\kph})\big)\big|,\,
\big|\lambda_d\big({\cal B}(\mU^-_{j,\kph})\big)\big|,\,\big|\lambda_d\big({\cal B}(\mU^+_{j,\kph})\big)\big|\Big\}.
\end{aligned}
\end{equation*}

As in the 1-D case, the global fluxes $\bm K^\pm_{\jph,k}=\bm K(\mU_{\jph,k}^\pm)$ and $\bm L^\pm_{j,\kph}=\bm L(\mU_{j,\kph}^\pm)$ in
\eref{5.2} are to be computed in such a way to ensure both WB and path-conservative properties of the resulting scheme. To this end, we
follow the approach from \cite{KLX_21} and describe how to compute the global fluxes $\bm K^\pm_{\jph,k}$ ($\bm L^\pm_{j,\kph}$ can be
computed in a similar manner). We first note that \eref{1.2b} admits steady-state solutions satisfying
\begin{equation*}
\bm F(\bm U)_x-B(\bm U)\bm U_x=\bm0, \quad {\rm and}\quad\bm G(\bm U)_y-C(\bm U)\bm U_y=\bm0,
\end{equation*}
which for many particular nonconservative systems can be written as
\begin{equation}
\bm F(\bm U)_x-B(\bm U)\bm U_x=M^x(\bm U)\bm E^x(\bm U)_x=\bm0,
\label{5.1b}
\end{equation}
where $M^x\in\mathbb R^{d\times d}$ and $\bm E^x$ is the vector of equilibrium variables. If $M^x$ is invertible (which is the case in most
practically interesting cases), the steady states satisfy the relation
\begin{equation*}
\bm E^x(\bm U)\equiv\widehat{\bm E}^x(y)
\end{equation*}
with $\widehat{\bm E}^x$ being a function of $y$ only. Taking this into account, we evaluate
\begin{equation*}
\bm K^\pm_\jph=\bm F^\pm_{\jph,k}-(\bm R^x)^\pm_{\jph,k},
\end{equation*}
where $\bm F^\pm_{\jph,k}:=\bm F(\bm U^\pm_{\jph,k})$, and the point values $\bm U^\pm_{\jph,k}$ are to be obtained using a 1-D
reconstruction of the equilibrium variables $\bm E^x$ to ensure the WB property. To this end, we first compute the point values
$\mE^x_{j,k}:=\mE^x(\mU_{j,k})$ and then apply a WENO-Z interpolation in the $x$-direction to evaluate $(\bm E^x)^\pm_{\jph,k}$. The
corresponding values $\bm U^\pm_{\jph,k}$ are then computed by solving the (nonlinear systems of) equations
\begin{equation*}
\bm E^x(\bm U^+_{\jph,k})=(\bm E^x)^+_{\jph,k}\quad\mbox{and}\quad\bm E^x(\bm U^-_{\jph,k})=(\bm E^x)^-_{\jph,k}
\end{equation*}
for $\bm U^+_{\jph,k}$ and $\bm U^-_{\jph,k}$, respectively. The point values of the global variables $(\bm R^x)^\pm_{\jph,k}$ in
\eref{5.1b} are then computed in a path-conservative way. We first select $\widehat x=x_\hf$, set $(\bm R^x)^-_{\hf,k}:=\bm0$, evaluate
\begin{equation}
(\bm R^x)^+_{\hf,k}=\bm B_{\bm\Psi,\hf,k},
\label{5.5}
\end{equation}
and then recursively obtain
\begin{equation}
(\bm R^x)^-_{\jph,k}=(\bm R^x)^+_{\jmh,k}+\bm B_{j,k},\quad(\bm R^x)^+_{\jph,k}=(\bm R^x)^-_{\jph,k}+\bm B_{\bm\Psi,\jph,k},
\label{5.6}
\end{equation}
for $j=1,\ldots,N_x$, $k=1,\ldots,N_y.$ In \eref{5.5} and \eref{5.6}, $\bm B_{j,k}$ and $\bm B_{\bm\Psi,\jph,k}$ are evaluated using
\eref{5.1b} and proper quadratures:
\begin{equation}
\begin{aligned}
&\bm B_{j,k}=\int\limits^{x_\jph}_{x_\jmh}B(\bm U)\bm U_x(x,y_k,t)\,{\rm d}x=\bm F_{\jph,k}^--\bm F_{\jmh,k}^+-
\int\limits^{x_\jph}_{x_\jmh}M^x(\bm U)\bm E^x(\bm U)_x(x,y_k,t)\,{\rm d}x,\\
&\bm B_{\bm\Psi,\jph,k}=\int\limits^1_0B\big(\bm\Psi_{\jph,k}(s)\big)\bm\Psi'_{\jph,k}(s)\,{\rm d}s\\
&\hspace{1.5cm}=\bm F_{\jph,k}^+-\bm F_{\jph,k}^--\hf\left[M^x\big(\bm U_{\jph,k}^+\big)+M^x\big(\bm U_{\jph,k}^-\big)\right]
\left[(\bm E^x)_{\jph,k}^+-(\bm E^x)_{\jph,k}^-\right],
\label{5.8}
\end{aligned}
\end{equation}
where $\bm\Psi_{\jph,k}(s):=\bm\Psi_{\jph,k}(s;\bm U^-_{\jph,k},\bm U^+_{\jph,k})$ is a sufficiently smooth path connecting the left- and
right-sided states $\bm U^-_{\jph,k}=\bm U\big(\bm E^-_{\jph,k}\big)$ and $\bm U^+_{\jph,k}=\bm U\big(\bm E^+_{\jph,k}\big)$. We refer the
readers to \cite{CKX23_WB} for more details on the derivation of \eref{5.8}.

In order to ensure the designed 2-D flux globalization based WB path-conservative A-WENO scheme is fifth-order accurate, the integral
$\bm B_{j,k}$ in \eref{5.8} needs to be evaluated using at least a fifth-order quadrature. In this paper, we use the Newton-Cotes quadrature
from \cite{Chu21}.

Finally, $(\bm K_{xx})_{\jph,k}$, $({\bm K_{xxxx}})_{\jph,k}$,  $(\bm L_{yy})_{j,\kph}$, and $({\bm L_{yyyy}})_{j,\kph}$ are approximations
of the second- and fourth-order spatial derivatives of $\bm K$ and $\bm L$ at $(x_\jph,y_k)$ and $(x_j,y_\kph)$, respectively, which can be
computed using the same FD approximations as in \eref{4.1}, namely, by
\begin{equation*}
\begin{aligned}
&(\mK_{xx})_{\jph,k}=\frac{1}{48(\dx)^2}\left(-5\mK_{j-2,k}+39\mK_{j-1,k}-34\mK_{j,k}-34\mK_{j+1,k}+39\mK_{j+2,k}-5\mK_{j+3,k}\right),\\
&(\mK_{xxxx})_{\jph,k}=\frac{1}{2(\dx)^4}\left(\mK_{j-2,k}-3\mK_{j-1,k}+2\mK_{j,k}+2\mK_{j+1,k}-3\mK_{j+2,k}+\mK_{j+3,k}\right), \\
&(\mL_{yy})_{j,\kph}=\frac{1}{48(\dy)^2}\left(-5\mL_{j,k-2}+39\mL_{j,k-1}-34\mL_{j,k}-34\mL_{j,k+1}+39\mL_{j,k+2}-5\mL_{j,k+3}\right),\\
&(\mL_{yyyy})_{j,\kph}=\frac{1}{2(\dy)^4}\left(\mL_{j,k-2}-3\mL_{j,k-1}+2\mL_{j,k}+2\mL_{j,k+1}-3\mL_{j,k+2}+\mL_{j,k+3}\right),
\end{aligned}
\end{equation*}
where $\mK_{j,k}=\mF(\mU_{j,k})-\bm R^x_{j,k}(\bm U)$, $\mL_{j,k}=\mG(\mU_{j,k})-\bm R^y_{j,k}(\bm U)$,  and
\begin{equation}
\begin{aligned}
&{\bm R^x_{j,k}}(\bm U)=\int\limits_{\widehat x}^{x_j}B({\bm U(x,y_k,t)})\bm U_x(x,y_k,t)\,{\rm d}x,\quad
{\bm R^y_{j,k}}(\bm U)=\int\limits_{\widehat y}^{y_k}C({\bm U(x_j,y_k,t)})\bm U_y(x_j,y,t)\,{\rm d}y.
\label{5.16}
\end{aligned}
\end{equation}
In order to obtain a fifth-order numerical scheme, the integrals in \eref{5.16} need to be evaluated within the fifth order of accuracy. We
refer the reader to, e.g., \cite{CKMZ23} for the fifth-order Newton-Cotes quadrature, which should be implemented in a recursive way.

\subsection{More Efficient Two-Dimensional Flux Globalization Based A-WENO Schemes}\label{sec5.2}
In this section, we extend a new, more efficient A-WENO scheme from \S\ref{sec4.2} to the nonconservative system \eref{1.2b}. We now
compute the high-order correction terms $(\bm K_{xx})_{\jph,k}$, $({\bm K_{xxxx}})_{\jph,k}$, $(\bm L_{yy})_{j,\kph}$, and
$({\bm L_{yyyy}})_{j,\kph}$ with the help of the same central difference formulae as in \S\ref{sec2.2}, \S\ref{sec3.2}, and \S\ref{sec4.2},
namely, by
\begin{equation}
\begin{aligned}
&(\mK_{xx})_{\jph,k}=\frac{1}{12(\dx)^2}\Big[-\bm{{\cal K}}^{\rm FV}_{j-\frac{3}{2},k}+16\bm{{\cal K}}^{\rm FV}_{\jmh,k}-
30\bm{{\cal K}}^{\rm FV}_{\jph,k}+16\bm{{\cal K}}^{\rm FV}_{j+\frac{3}{2},k}-\bm{{\cal K}}^{\rm FV}_{j+\frac{5}{2},k}\Big],\\
&(\mK_{xxxx})_{\jph,k}=\frac{1}{(\dx)^4}\Big[\bm{{\cal K}}^{\rm FV}_{j-\frac{3}{2},k}-4\bm{{\cal K}}^{\rm FV}_{\jmh,k}+
6\bm{{\cal K}}^{\rm FV}_{\jph,k}-4\bm{{\cal K}}^{\rm FV}_{j+\frac{3}{2},k}+\bm{{\cal K}}^{\rm FV}_{j+\frac{5}{2},k}\Big],\\
&(\mL_{yy})_{j,\kph}=\frac{1}{12(\dy)^2}\Big[-\bm{{\cal L}}^{\rm FV}_{j,k-\frac{3}{2}}+16\bm{{\cal L}}^{\rm FV}_{j,\kmh}-
30\bm{{\cal L}}^{\rm FV}_{j,\kph}+16\bm{{\cal L}}^{\rm FV}_{j,k+\frac{3}{2}}-\bm{{\cal L}}^{\rm FV}_{j,k+\frac{5}{2}}\Big],\\
&(\mL_{yyyy})_{j,\kph}=\frac{1}{(\dy)^4}\Big[\bm{{\cal L}}^{\rm FV}_{j,k-\frac{3}{2}}-4\bm{{\cal L}}^{\rm FV}_{j,\kmh}+
6\bm{{\cal L}}^{\rm FV}_{j,\kph}-4\bm{{\cal L}}^{\rm FV}_{j,k+\frac{3}{2}}+\bm{{\cal L}}^{\rm FV}_{j,k+\frac{5}{2}}\Big].
\end{aligned}
\label{5.8f}
\end{equation}
\begin{remark}
As in the conservative case, the presented flux globalization based WB path-conservative A-WENO scheme can be straightforwardly extended to
higher than the fifth order of accuracy via higher-order Taylor expansions.
\end{remark}

\section{Numerical Examples}\label{sec6}
For the sake of brevity, the A-WENO schemes from \S\ref{sec2.1}, \S\ref{sec3.1}, \S\ref{sec4.1}, and \S\ref{sec5.1} will be
referred to as the Old Schemes, while the A-WENO schemes from \S\ref{sec2.2}, \S\ref{sec3.2}, \S\ref{sec4.2}, and
\S\ref{sec5.2} will be referred to as the New Schemes. In this section, we will test the New and Old Schemes on several numerical examples.

We numerically integrate the ODE systems \eref{2.1a}, \eref{3.1}, \eref{4.1a}, and \eref{5.1} by the three-stage third-order
strong stability preserving (SSP) Runge-Kutta method (see, e.g., \cite{Gottlieb11,Gottlieb12}) and use the CFL number 0.45,
that is, the time step $\dt$ is selected adaptively as $\dt=\frac{0.45\dx}{\max\limits_ja_\jph}$ in the 1-D examples and
$\dt=0.45\max\big\{\frac{\dx}{\max\limits_{j,k}a_{\jph,k}},\frac{\dy}{\max\limits_{j,k}a_{j,\kph}}\big\}$ in the 2-D examples. In Example
4, where we test the accuracy of the studied A-WENO schemes, we take
$\dt=\dfrac{0.45(\dx)^{\frac{5}{3}}}{\max\limits_ja_\jph}$ to ensure that the time errors do not dominate the spatial ones.

In all of the numerical examples, boundary conditions are implemented using a ghost point technique. We stress that both the
point values of the computed solution and numerical fluxes used in the computation of the high-order correction terms in \eref{2.3},
\eref{3.10f}, \eref{4.4f}, and \eref{5.8f} are extrapolated across the boundaries. Alternatively, one can expand the stencil and then to
extrapolate the point values only. This approach can be favorable in the case when the the fluxes of the system at hand are very
complicated.

All of the numerical experiments were conducted on a workstation equipped with an Intel(R) Core(TM) i7-9750H CPU at 2.6 GHz and 32 GB of
RAM. The simulations were implemented in Fortran using GCC version 14.2.0 compiler suite. The reported CPU times were averaged over thirty
independent runs to ensure reproducibility and minimize variability due to system processes.

\subsection{One-Dimensional Examples}
In this section, we present several 1-D examples to compare the results computed by the New and Old Schemes, along with the CPU times
consumed by each scheme.

\subsubsection{Scalar Equations}
\paragraph*{Example 1---Burgers Equation.} In the first example, we consider the inviscid Burgers equation
$$
u_t+\Big(\dfrac{u^2}{2}\Big)_x=0
$$
subject to the 1-periodic initial conditions
$$
u(x,0)=\frac{1}{4}+\hf\sin(2\pi x).
$$

\begin{figure}[ht!]
\centerline{
\includegraphics[trim=0.8cm 0.4cm 1.3cm 0.6cm, clip, width=6.0cm]{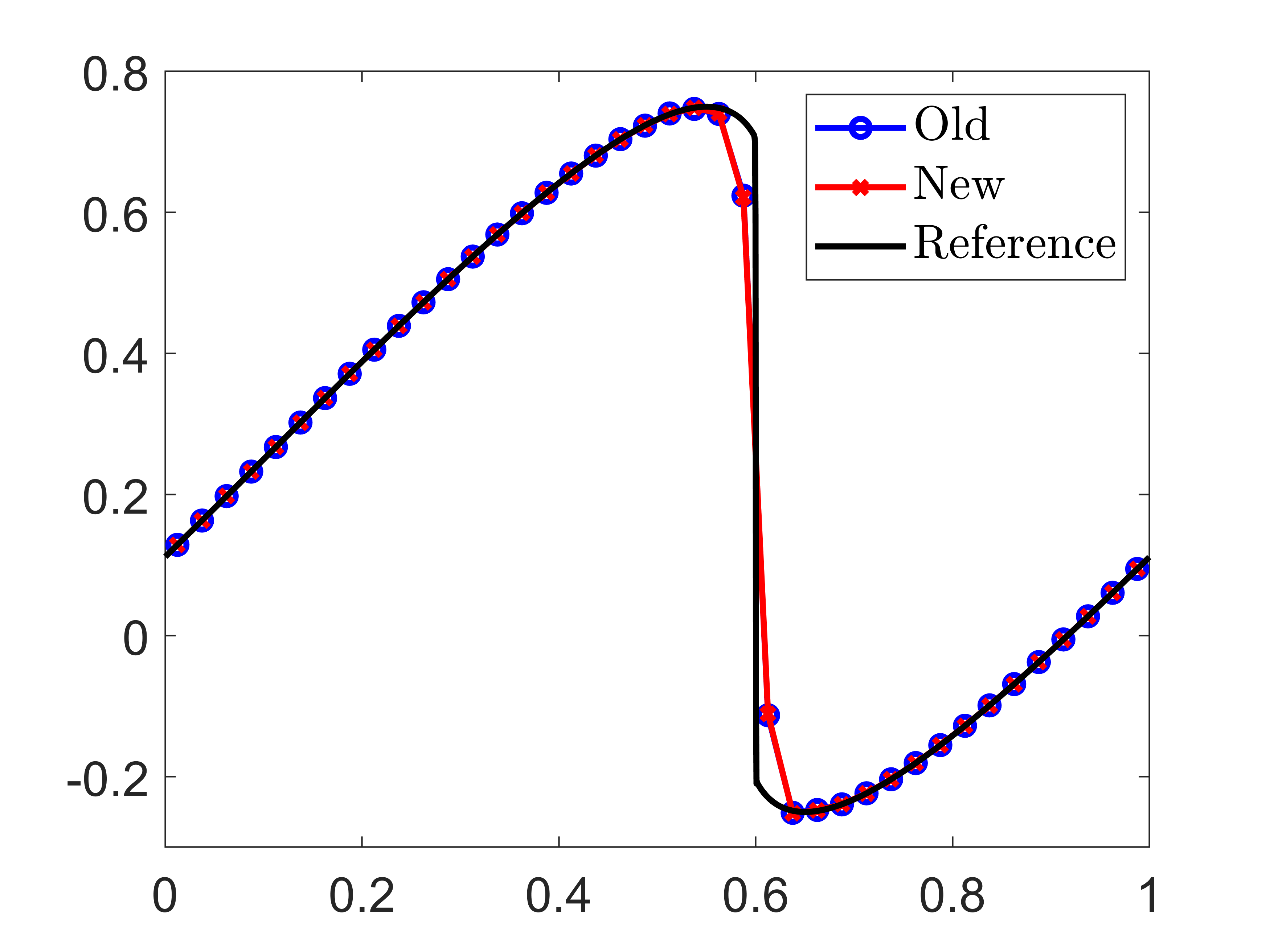}}
\caption{\sf Example 1: Numerical solutions computed by the New and Old Schemes.\label{fig41}}
\end{figure}

We first compute the numerical solutions by both the New and Old Schemes on the uniform mesh with $\dx=1/40$ in the computational domain
$[0,1]$ until the final time $t=0.4$. The obtained numerical results are plotted in Figure \ref{fig41} along with the reference solution
computed by the Old Scheme on a much finer mesh with $\dx=1/2000$. As one can clearly see, the results obtained by the New and Old Schemes
are almost identical.

We then measure the CPU times consumed by both of the studied schemes on a much finer mesh with $\dx=1/40000$ and present the
obtained CPU times in Table \ref{tab_ex1}. As one can see, the CPU time consumed by the Old Scheme is about $4.4\%$ larger than the CPU
time consumed by the New Scheme. The difference is relatively small since the flux function $f(u)=u^2/2$ is very simple in this example.
\begin{table}[ht!]
\centering
\begin{tabular}{|c|c|c|}
\hline
Old Scheme&New Scheme&Ratio\\
\hline
68.82&65.94&1.044\\
\hline
\end{tabular}
\caption{\sf Example 1: The CPU times (in seconds) consumed by the studied schemes.\label{tab_ex1}}
\end{table}

\begin{figure}[ht!]
\centerline{
\includegraphics[trim=0.8cm 0.4cm 1.3cm 0.6cm, clip, width=6.0cm]{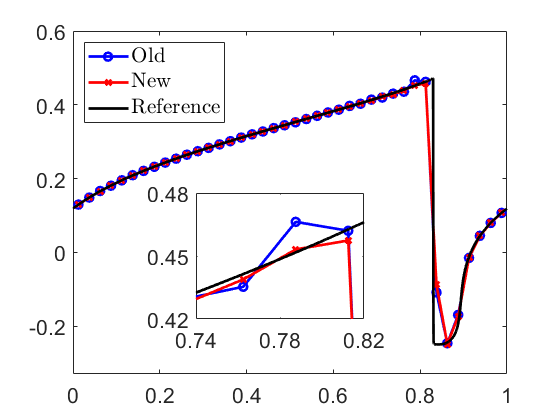}}
\caption{\sf Example 2: Numerical solutions computed by the New and Old Schemes with the zoom at the left part of the shock area.
\label{fig42}}
\end{figure}

\paragraph*{Example 2---Buckley-Leverett Equation with Gravitational Effects.}
In the second example, we consider the Buckley-Leverett equation with gravitational effects, which reads as
$$
u_t+\bigg[\frac{u^2}{u^2+(1-u)^2}(1-k(1-u)^2)\bigg]_x=0.
$$
We take $k=1$, consider the 1-periodic initial conditions
$$
u(x,0)=\frac{1}{4}+\hf\sin(2\pi x),
$$
and compute the numerical solutions by both the New and Old Schemes on the uniform mesh with $\dx=1/40$ in the computational domain $[0,1]$
until the final time $t=0.4$. In Figure \ref{fig42}, we present the obtained results along with the reference solution computed by the Old
Scheme on a much finer mesh with $\dx=1/2000$. As in Example 1, the New and Old solutions almost coincide. However, when we zoom at the area
near the shock, one can see that the New Scheme produces less oscillatory solution than the Old Scheme.

Next, we measure the CPU times consumed by both of the studied schemes on a much finer mesh with $\dx=1/20000$. The obtained
results are reported in Table \ref{tab_ex2}, where one can see that the CPU time consumed by the Old Scheme is about $7.6\%$ larger than
the CPU time consumed by the New Scheme. The efficiency gain is now bigger than in Example 1. This suggests that one can save more CPU time
when the flux function is more complex.
\begin{table}[ht!]
\centering
\begin{tabular}{|c|c|c|}
\hline
Old Scheme&New Scheme&Ratio\\
\hline
111.97&104.06&1.076\\
\hline
\end{tabular}
\caption{\sf Example 2: The CPU times (in seconds) consumed by the studied schemes.\label{tab_ex2}}
\end{table}

\paragraph*{Example 3---Scalar Equation with a Source Term.}
In the third example taken from \cite{NK06}, we consider the following nonconservative scalar equation:
\begin{equation}
u_t+f(u)_x+z_xu=0,\quad f(u)=\frac{u^2}{2},\quad z(x)=\begin{cases}-\cos(\pi x)&\mbox{if}~x\in(\frac{3}{2},\frac{5}{2}),\\
0&\mbox{otherwise}.
\end{cases}
\label{equ4.1}
\end{equation}
We take the constant initial datum $u(x,0)\equiv1$ and impose the Dirichlet boundary conditions $u(0,t)=2$ and $u(4,t)=1$ on the interval
$[0,4]$. It is easy to show that at the steady states $E:=u+z\equiv{\rm Const}$ and a scalar version of equation \eref{3.1b},
$$
f(u)_x+z_xu=M(u)E_x=0,
$$
is satisfied with $M(u)=u$.

\begin{figure}[ht!]
\centerline{\includegraphics[trim=0.9cm 0.4cm 1.3cm 0.2cm, clip, width=6.0cm]{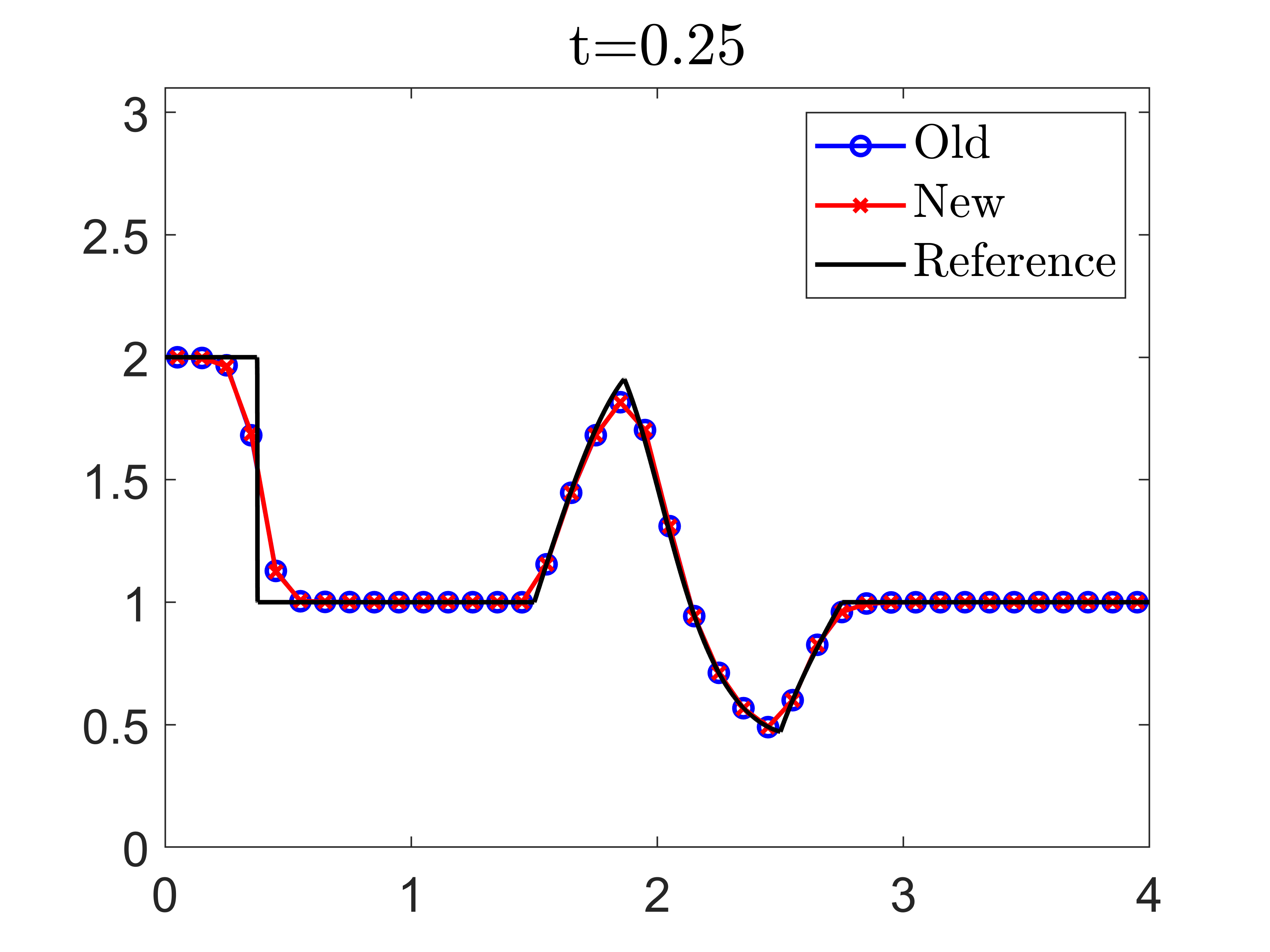}\hspace*{1cm}
            \includegraphics[trim=0.9cm 0.4cm 1.3cm 0.2cm, clip, width=6.0cm]{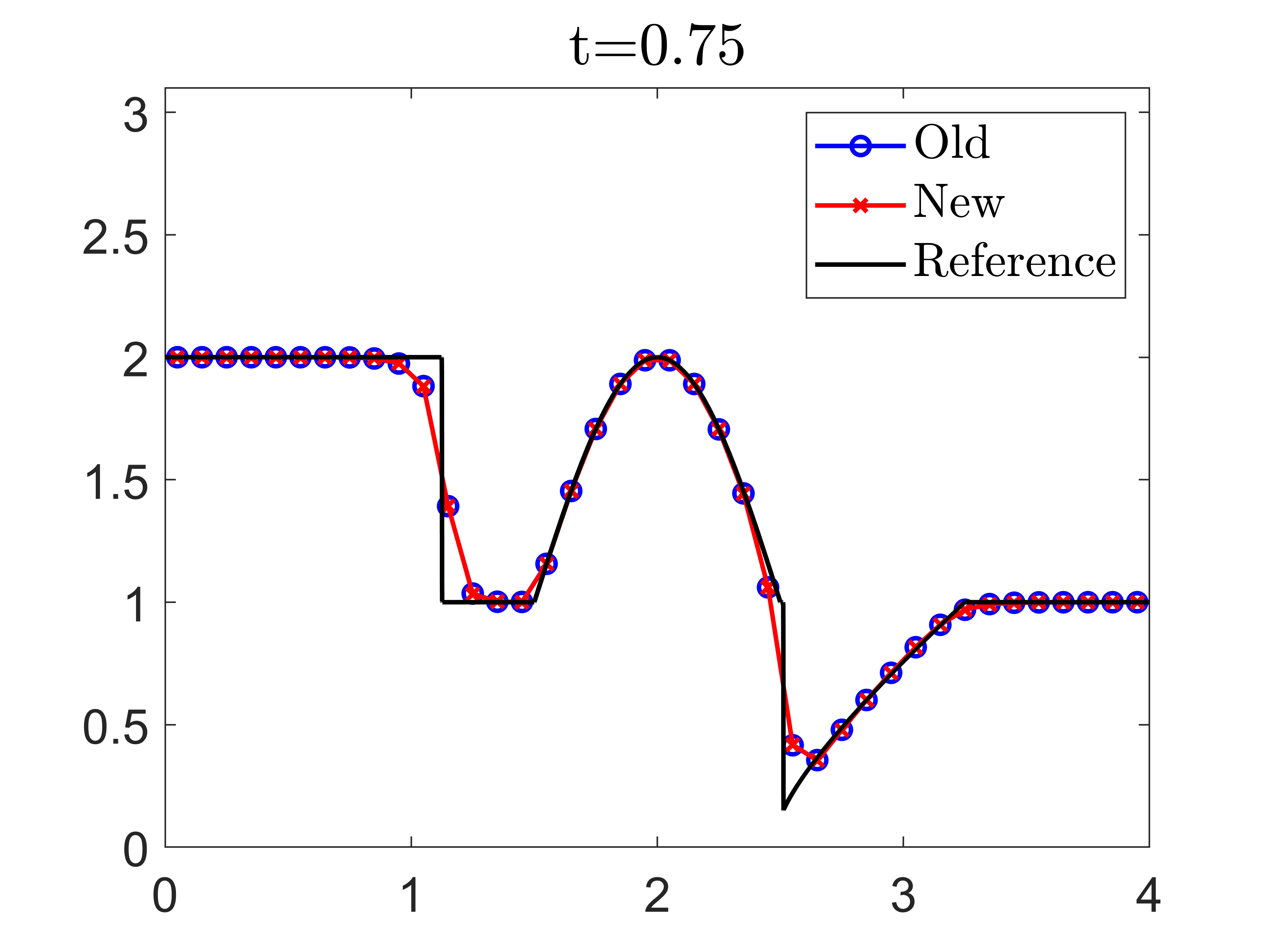}}
\vskip10pt
\centerline{\includegraphics[trim=0.9cm 0.4cm 1.3cm 0.2cm, clip, width=6.0cm]{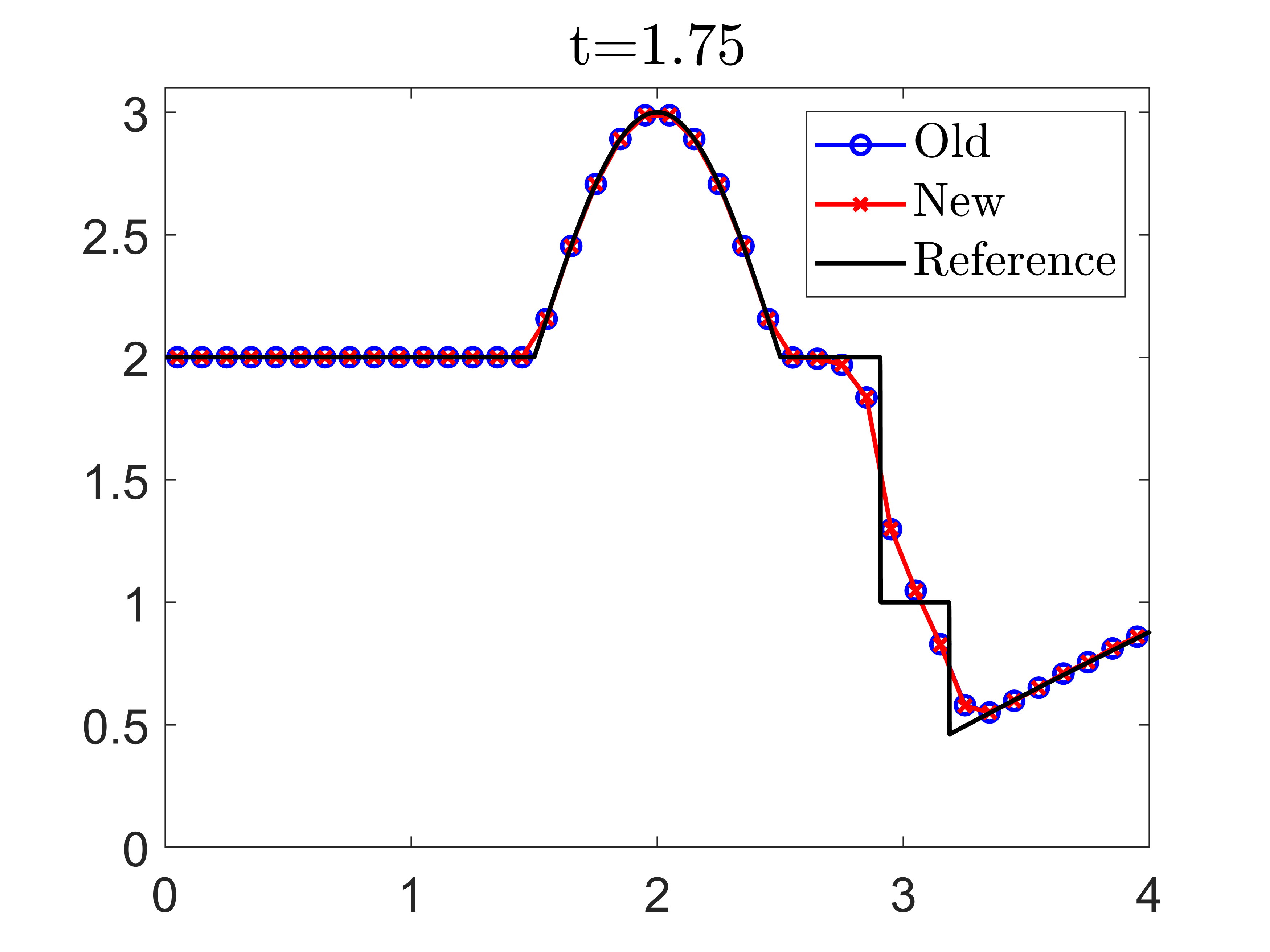}\hspace*{1cm}
            \includegraphics[trim=0.9cm 0.4cm 1.3cm 0.2cm, clip, width=6.0cm]{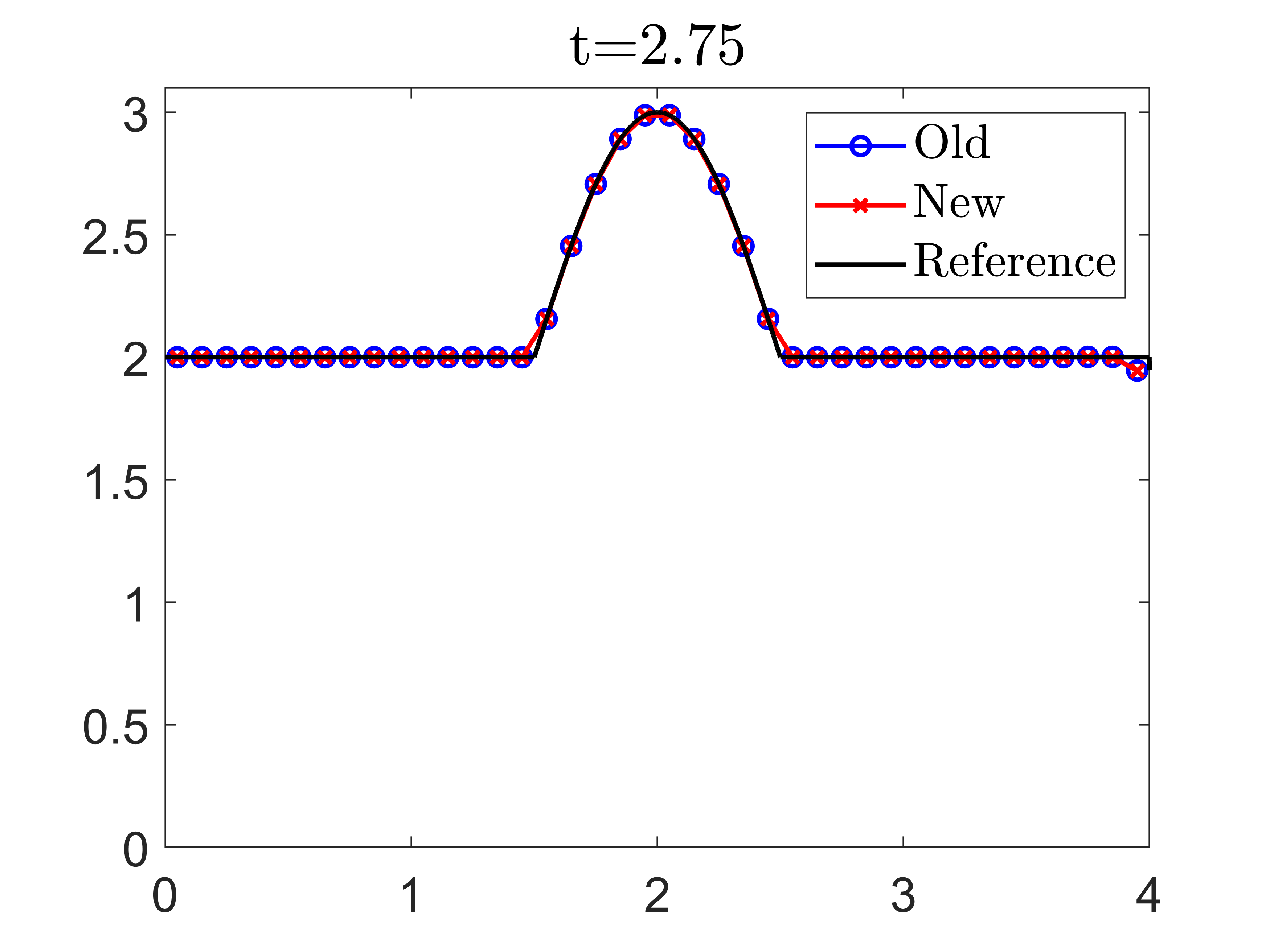}}
\caption{\sf Example 3: Numerical solutions computed by the New and Old Schemes at different times.\label{fig43}}
\end{figure}

We follow \S\ref{sec3} and design the Old and New Schemes for \eref{equ4.1}. We then compute the numerical solutions by both the obtained
schemes on the uniform mesh with $\dx=1/10$ until the final time $t=2.75$ and present the snapshots of the computed solutions at $t=0.25$,
0.75, 1.75, and 2.75 in Figure \ref{fig43} along with the reference solution computed by the Old Scheme on a much finer mesh with
$\dx=1/10000$. Compared with the results reported in \cite[\S5]{NK06}, our solutions are sharper (as expected since our schemes are higher
order). At the same time, one can clearly see that the results obtained by the New and Old Schemes are almost identical.

We then measure the CPU times consumed by both of the studied schemes on a much finer mesh with $\dx=1/1500$ and present the
obtained results in Table \ref{tab_ex3}, which shows that the CPU times consumed by the Old Scheme are about $15.1\%$ larger than the CPU
time consumed by the New Scheme. The difference is larger than in Examples 1 and 2. This suggests that one can save much more CPU time by
applying the new A-WENO schemes to nonconservative hyperbolic PDEs.
\begin{table}[ht!]
\centering
\begin{tabular}{|c|c|c|}
\hline
Old Scheme&New Scheme&Ratio\\
\hline
139.85&121.54&1.151\\
\hline 
\end{tabular}
\caption{\sf Example 3: The CPU times (in seconds) consumed by the studied schemes.\label{tab_ex3}}
\end{table}

\subsubsection{One-Dimensional Euler Equations of Gas Dynamics}
In this section, we consider the 1-D Euler equations of gas dynamics:
\begin{equation}
\begin{aligned}
&\rho_t+(\rho u)_x=0,\\
&(\rho u)_t+(\rho u^2+p)_x=0,\\
&E_t+\left[u(E+p)\right]_x=0,
\end{aligned}
\label{4.2}
\end{equation}
where $\rho$, $u$, $p$, and $E$ are the density, velocity, pressure and total energy, respectively. The system is completed through the
equation of state for the ideal gas:
\begin{equation}
p=(\gamma-1)\Big[E-\hf\rho u^2\Big],
\label{4.5}
\end{equation}
where the parameter $\gamma$ represents the specific heat ratio (we take $\gamma=1.4)$. When applying the studied A-WENO schemes to
\eref{4.2}--\eref{4.5}, we perform the WENO-Z interpolation in the local characteristic variables using the local characteristic
decomposition, whose detailed description can be found in \cite{CCK23_Adaptive}.

In Example 4, we test the accuracy of the New and Old Schemes and demonstrate that both the errors and experimental convergence rates are
almost the same. In Examples 5 and 6, we consider more challenging numerical examples and demonstrate that the New Scheme is capable of
achieving the same resolution as the Old Scheme.

\paragraph*{Example 4---Accuracy Test.}
In this example taken from \cite{KKOKC}, we consider the following smooth initial data:
\begin{equation*}
u(x,0)=\sin\Big(\frac{\pi x}{5}+\frac{\pi}{4}\Big),\quad\rho(x,0)=
\bigg[\frac{\gamma-1}{2\sqrt{\gamma}}\left(u(x,0)+10\right)\bigg]^{\frac{2}{\gamma-1}},\quad p(x,0)=\rho^\gamma(x,0),
\end{equation*}
subject to the 10-periodic boundary conditions in the computational domain $[0,10]$.

\begin{figure}[ht!]
\centerline{
\includegraphics[trim=0.8cm 0.4cm 1.1cm 0.6cm, clip, width=6.0cm]{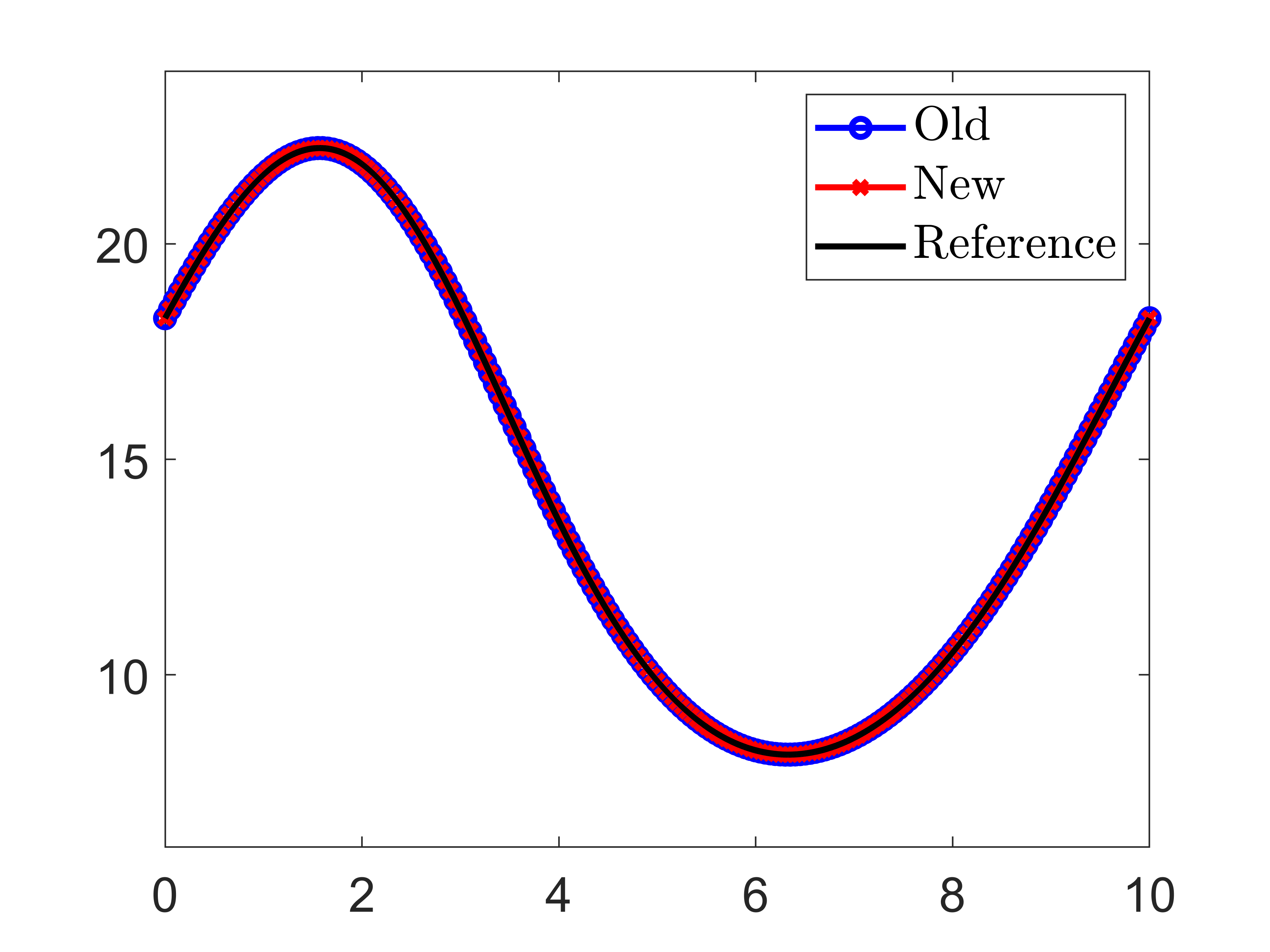}}
\caption{\sf Example 4: Density $\rho$ computed by the New and Old Schemes with $\dx=1/20$.\label{fig44}}
\end{figure}

We compute the numerical solutions by both the New and Old Schemes until the final time $t=0.1$ using both the New and Old Schemes on a
sequence of uniform meshes with $\dx=1/20$, $1/40$, $1/80$, $1/160$, $1/320$, and $1/640$. The densities computed with $\dx=1/20$ are
plotted in Figure \ref{fig44} along with the reference solution computed by the Old Scheme on a much finer mesh with $\dx=1/320$, where one
can clearly see that the results obtained by the New and Old Schemes are almost identical.

We then compute the $L^1$-errors and estimate the experimental convergence rates using the following Runge formulae, which are based on the
solutions computed on the three consecutive uniform grids with the mesh sizes $\dx$, $2\dx$, and $4\dx$ and denoted by $(\cdot)^{\dx}$,
$(\cdot)^{2\dx}$, and $(\cdot)^{4\dx}$, respectively:
$$
{\rm Error}(\dx)\approx\frac{\delta_{12}^2}{|\delta_{12}-\delta_{24}|},\quad
{\rm Rate}(\dx)\approx\log_2\left(\frac{\delta_{24}}{\delta_{12}}\right).
$$
Here, $\delta_{12}:=\|(\cdot)^{\dx}-(\cdot)^{2\dx}\|_{L^1}$ and $\delta_{24}:=\|(\cdot)^{2\dx}-(\cdot)^{4\dx}\|_{L^1}$. The obtained results
for the density and total energy are reported in Table \ref{tab1}, where one can clearly see that the fifth order of accuracy is achieved by
the New Scheme and the errors are the same as those in the Old Scheme results.
\begin{table}[ht!]
\centering
\begin{tabular}{|c|cc|cc|cc|cc|}
\hline
\multirow{3}{1em}{$\dx$}&\multicolumn{4}{c|}{New A-WENO Scheme}&\multicolumn{4}{c|}{Old A-WENO Scheme}\\
\cline{2-9}
&\multicolumn{2}{c|}{$\rho$}&\multicolumn{2}{c|}{$E$}&\multicolumn{2}{c|}{$\rho$}&\multicolumn{2}{c|}{$E$}\\
\cline{2-9}&Error&Rate&Error&Rate&Error&Rate&Error&Rate\\
\hline
$1/160$&1.44e-09&4.79&2.13e-08&4.76&1.44e-09&4.79&2.13e-08&4.76\\
$1/320$&3.88e-11&5.00&5.65e-10&4.99&3.88e-11&5.00&5.65e-10&4.99\\
$1/640$&1.25e-12&4.98&1.81e-11&4.98&1.25e-12&4.98&1.81e-11&4.98\\
\hline
\end{tabular}
\caption{\sf Example 4: The $L^1$-errors and experimental convergence rates.\label{tab1}}
\end{table}

\begin{figure}[ht!]
\centerline{
\includegraphics[trim=0.7cm 0.3cm 1.2cm 0.6cm, clip, width=6.cm]{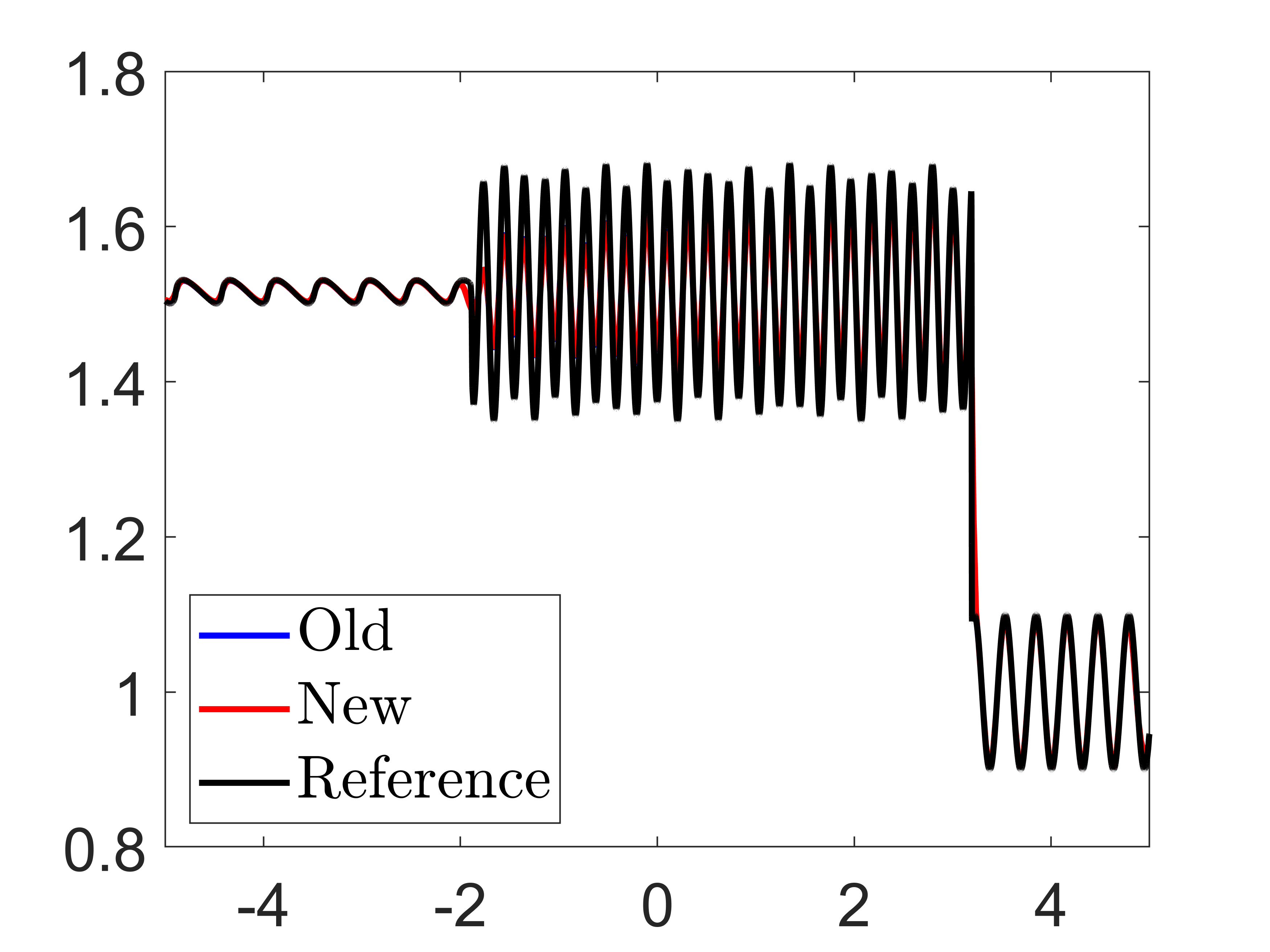}\hspace*{1.0cm}
\includegraphics[trim=0.7cm 0.3cm 1.2cm 0.6cm, clip, width=6.cm]{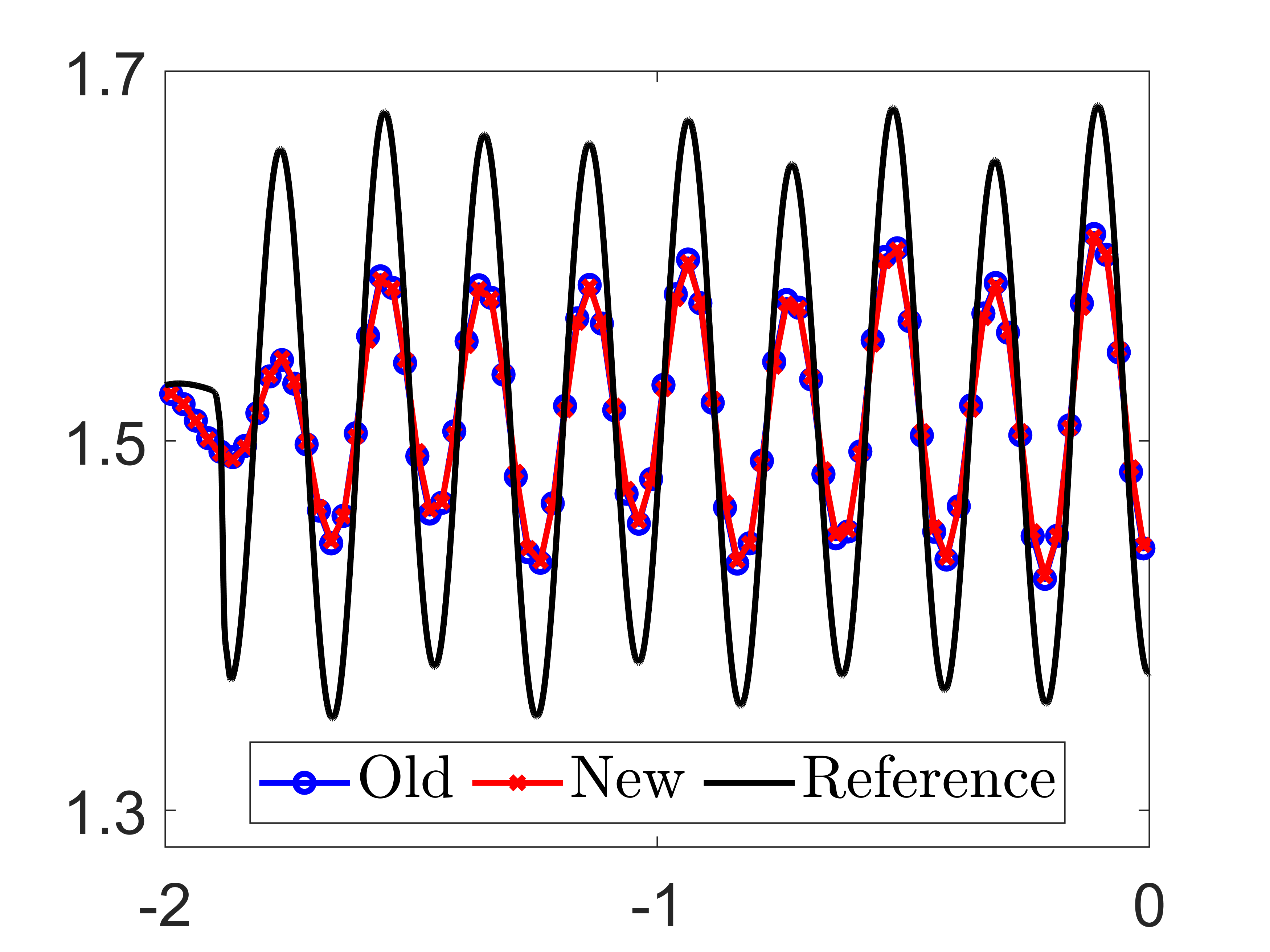}}
\caption{\sf Example 5: Density computed by the New and Old Schemes (left) and zoom at $x\in[-2,0]$ (right).\label{fig2}}
\end{figure}

\paragraph*{Example 5---Shock Entropy Problem.}
In this example, we consider the shock-entropy problem introduced in \cite{Shu88}. The initial data are
\begin{equation*}
(\rho,u,p)(x,0)=\begin{cases}
(1.51695,0.523346,1.805),&x<-4.5,\\
(1+0.1\sin(20x),0,1),&x>-4.5.
\end{cases}
\end{equation*}
We impose free boundary condition at the both ends of the computational domain $[-5,5]$ and compute the numerical solutions by both the New
and Old Schemes until the final time $t=5$ on a uniform mesh with $\dx=1/40$. The obtained numerical results are plotted in Figure
\ref{fig2} along with the reference solution computed by the Old Scheme on a much finer mesh with $\dx=1/800$. As one can see from Figure
\ref{fig2} (right), where we zoom at the area where the solution has smooth oscillatory structures, the results obtained by the New and Old
Schemes almost coincide, which suggests that the New Scheme is as accurate as the Old one.

We measure the CPU times consumed by both of the studied schemes on a much finer mesh with $\dx=1/600$, and report the obtained results in
Table \ref{tab_ex5},  where one can see that the CPU time consumed by the Old Scheme is about $2.2\%$ larger than the CPU time consumed by
the New Scheme. The efficiency gain is now smaller since we have used the local characteristic decomposition in WENO-Z reconstructions
(see, e.g., \cite{CCHKL_22,Don20,JSZ,Qiu02,Shu20,wang18} and references therein), which makes the computational cost of the A-WENO
correction terms relatively smaller.
\begin{table}[ht!]
\centering
\begin{tabular}{|c|c|c|}
\hline
Old Scheme&New Scheme&Ratio\\
\hline
127.18&124.44&1.022\\
\hline 
\end{tabular}
\caption{\sf Example 5: The CPU times (in seconds) consumed by the studied schemes.\label{tab_ex5}}
\end{table}

\paragraph*{Example 6---Blast Wave Problem.}
In this example taken from \cite{Woodward88}, we consider the following initial conditions:
\begin{equation*}
(\rho, u,p)(x,0)=\begin{cases}
(1,0,1000),&x<0.1,\\
(1,0,0.01),&0.1\le x\le0.9,\\
(1,0,100),&x>0.9,
\end{cases}
\end{equation*}
subject to the solid wall boundary conditions imposed at the both ends of the computational domain $[0,1]$.

\begin{figure}[ht!]
\centerline{
\includegraphics[trim=1.3cm 0.2cm 1.3cm 0.8cm, clip, width=6.cm]{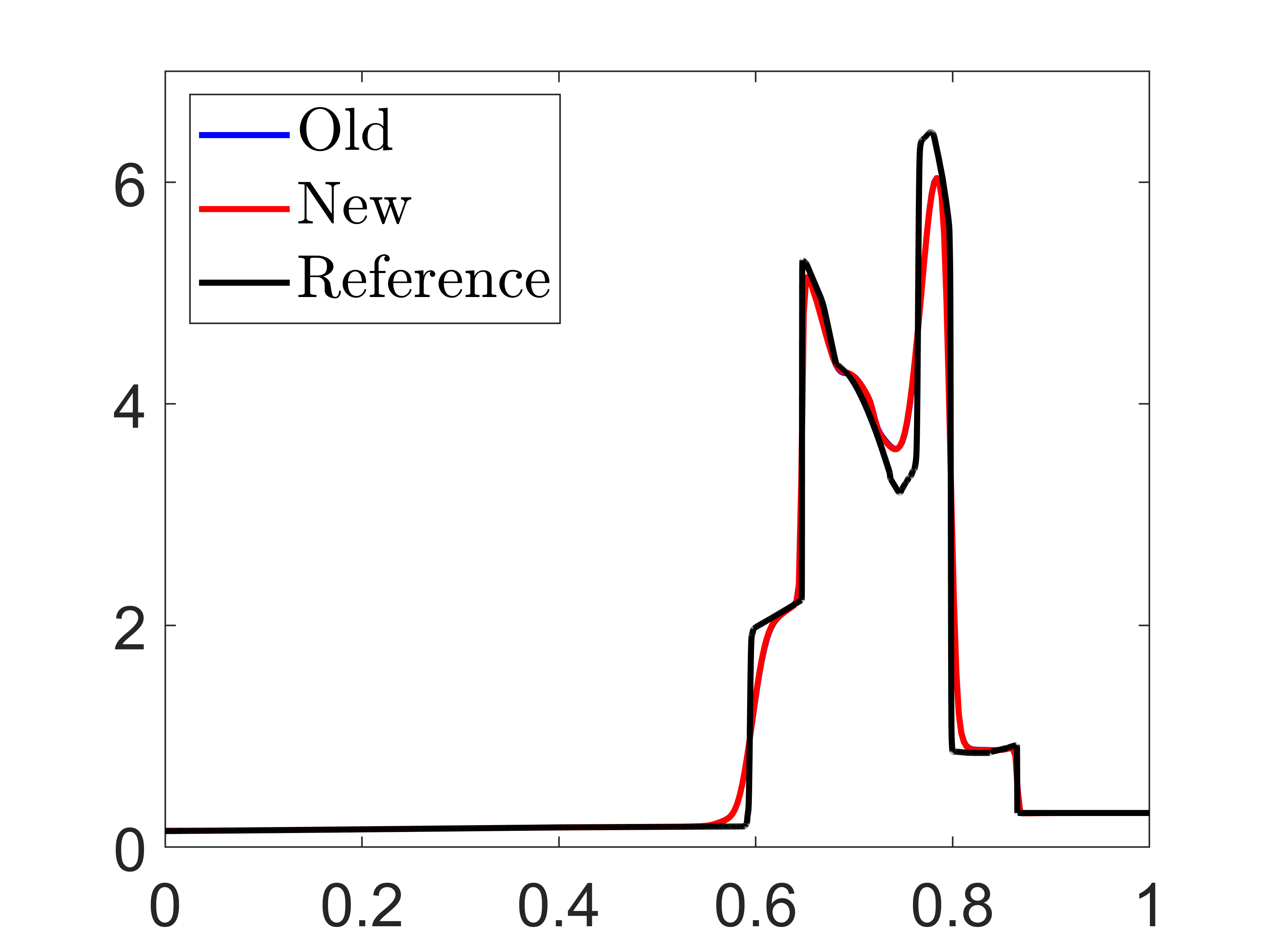}\hspace*{1.0cm}
\includegraphics[trim=1.3cm 0.2cm 1.3cm 0.8cm, clip, width=6.cm]{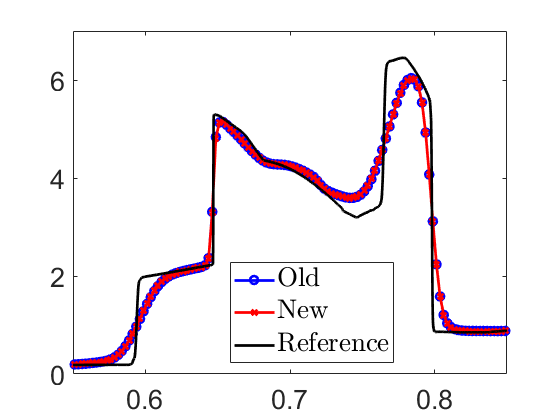}}
\caption{\sf Example 6: Density computed by the New and Old Schemes (left) and zoom at $x\in[0.55,0.85]$ (right).\label{fig3}}
\end{figure}

We compute the numerical solutions until the final time $t=0.038$ by both the New and Old Schemes on a uniform mesh with $\dx=1/400$ and
plot the obtained results in Figure \ref{fig3} along with the reference solution computed by the Old Scheme on a much finer grid with
$\dx=1/8000$. Once again, one can see that the New and Old Schemes achieve the same resolution.

We then measure the CPU times consumed by both of the studied schemes on a much finer mesh with $\dx=1/4000$. The obtained results are
reported in Table \ref{tab_ex6}, where one can see that the CPU time consumed by the Old Scheme is about $2.3\%$ larger than the CPU time
consumed by the New Scheme. 
\begin{table}[ht!]
\centering
\begin{tabular}{|c|c|c|}
\hline
Old Scheme&New Scheme&Ratio\\
\hline
93.92&91.78&1.023\\
\hline 
\end{tabular}
\caption{\sf Example 6: The CPU times (in seconds) consumed by the studied schemes.\label{tab_ex6}}
\end{table}

\subsection{Two-Dimensional Examples}
In this section, we present two 2-D examples to compare the results computed by the New and Old Schemes, as well as to evaluate the CPU
times consumed by the studied schemes.

\subsubsection{Two-Dimensional Euler Equations of Gas Dynamics}
In this section, we consider the 2-D Euler equations of gas dynamics, which read as
\begin{equation}
\begin{aligned}
&\rho_t+(\rho u)_x+(\rho v)_y=0,\\
&(\rho u)_t+(\rho u^2 +p)_x+(\rho uv)_y=0,\\
&(\rho v)_t+(\rho uv)_x+(\rho v^2+p)_y=0,\\
&E_t+\left[u(E+p)\right]_x+\left[v(E+p)\right]_y=0,
\end{aligned}
\label{6.6}
\end{equation}
where $v$ is the $y$-component of the velocity, and the rest of the notations are the same as in the 1-D case. The system is completed
through the following equations of state for the ideal gas:
\begin{equation*}
p=(\gamma-1)\Big[E-\frac{\rho}{2}(u^2+v^2)\Big].
\end{equation*}

\paragraph*{Example 7---Implosion Problem.}
In this example taken from \cite{Liska03}, we consider the implosion problem with the following initial conditions:
\begin{equation*}
(\rho,u,v,p)(x,y,0)=\begin{cases}
(0.125,0,0,0.14),&|x|+|y|<0.15,\\
(1,0,0,1),&\mbox{otherwise},
\end{cases}
\end{equation*}
prescribed in $[0,0.3]\times[0,0.3]$ subject to the solid wall boundary conditions at $x=0$ and $y=0$ and free boundary conditions at
$x=0.3$ and $y=0.3$.

\begin{figure}[ht!]
\centerline{\includegraphics[trim=5.0cm 3.4cm 1.9cm 2.5cm, clip, width=14cm]{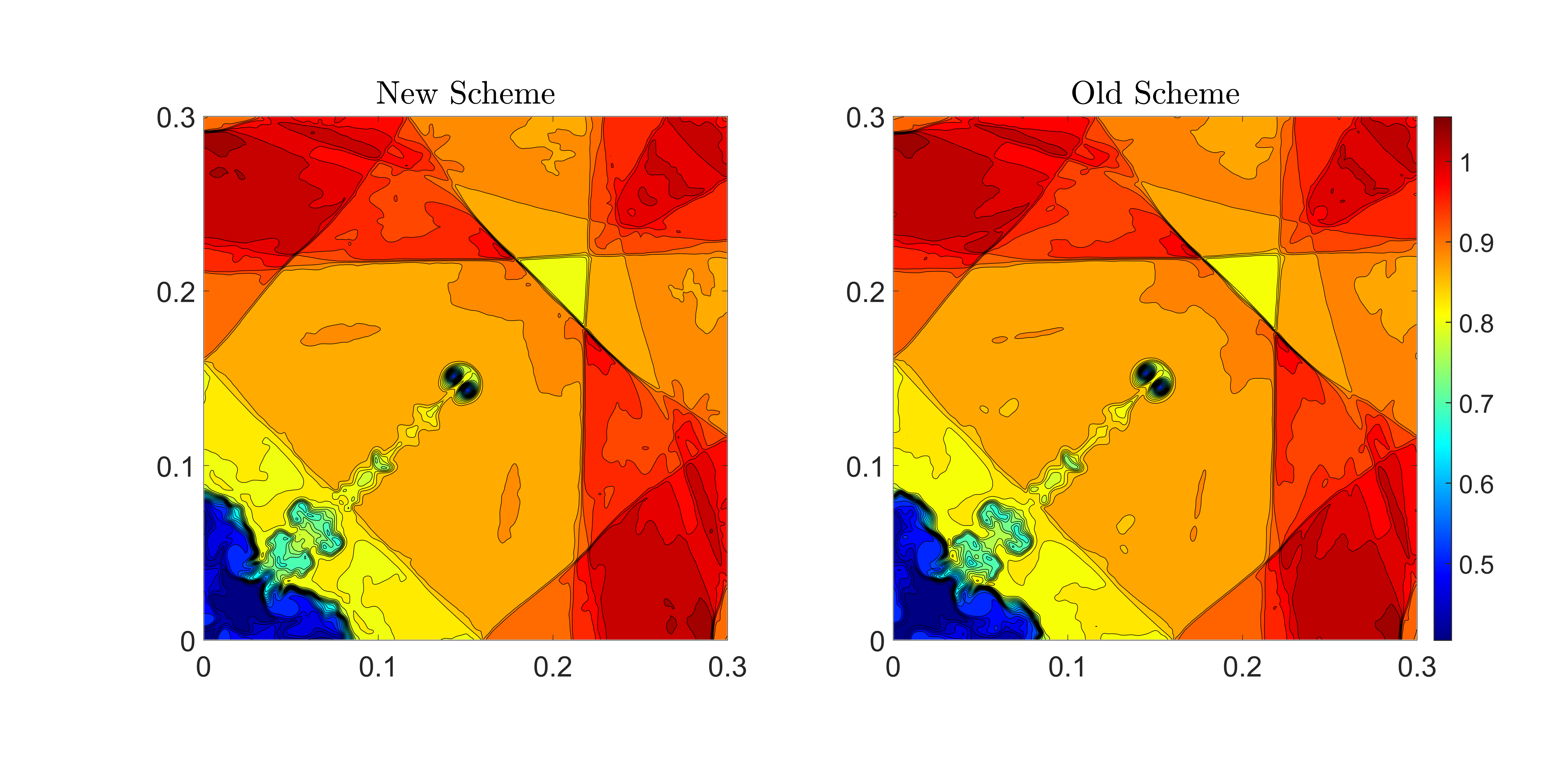}}
\caption{\sf Example 7: Density computed by the New (left) and Old (right) Schemes.\label{fig12}}
\end{figure}

We compute the numerical solutions by the New and Old Schemes until the final time $t=2.5$ on the uniform mesh with $\dx=\dy=3/4000$ and
plot the obtained results in Figure \ref{fig12}. As observed, the jets generated by the New and Old schemes propagate almost identically in
the direction of $y=x$, clearly indicating that the New and Old Schemes achieve the same resolution. 

Next, we measure the CPU times consumed by the New and Old Schemes. The obtained results are reported in Table \ref{tab_ex7}, where one can
see that the CPU time consumed by the Old Scheme is about $2.0\%$ larger than the CPU time consumed by the New Scheme. 
\begin{table}[ht!]
\centering
\begin{tabular}{|c|c|c|}
\hline
Old Scheme&New Scheme&Ratio\\
\hline
11690.44&11462.44&1.020\\
\hline
\end{tabular}
\caption{\sf Example 7: The CPU times (in seconds) consumed by the studied schemes.\label{tab_ex7}}
\end{table}

\subsubsection{Two-Dimensional $\gamma$-Based Multifluid System}
In this section, we consider the 2-D Euler equations of gas dynamics \eref{6.6}, supplemented with the following two equations for
$\Gamma:=1/(\gamma-1)$ and $\Pi:=\gamma\pi_\infty/(\gamma-1)$, which are used to track the interfaces between different components of the
multifluid:
\begin{equation*}
\begin{aligned}
&\Gamma_t+(u\Gamma)_x+(v\Gamma)_y=\Gamma(u_x+v_y),\\
&\Pi_t+(u\Pi)_x+(v\Pi)_y=\Pi(u_x+v_y).
\end{aligned}
\end{equation*}
The system is completed through the following stiff equations of state:
\begin{equation*}
p=(\gamma-1)\Big[E-\frac{\rho}{2}(u^2+v^2)\Big]-\gamma\pi_\infty,
\end{equation*}
where $\pi_\infty$ is a stiffness parameter.

\paragraph*{Example 8---A Cylindrical Explosion Problem.} In the last 2-D example taken from \cite{CKX24,CHK25,XL17}, we consider the case
where a cylindrical explosive source is located between an air-water interface and an impermeable wall. The initial conditions are given by
\begin{equation*}
(\rho,u,v,p;\gamma,\pi_\infty)\Big|_{(x,y,0)}=\begin{cases}(1.27,0,0,8290;2,0),&\mbox{$(x-5)^2+(y-2)^2<1$},\\
(0.02,0,0,1;1.4,0),&\mbox{$y>4$},\\(1,0,0,1;7.15,3309),&\mbox{otherwise},\end{cases}
\end{equation*}
the solid wall boundary conditions are imposed at the bottom, and the free boundary conditions are prescribed on the other sides of the
computational domain $[0,10]\times[0,6]$.

\begin{figure}[ht!]
\centerline{\includegraphics[trim=1.2cm 1.5cm 1.1cm 1.1cm, clip, width=7cm]{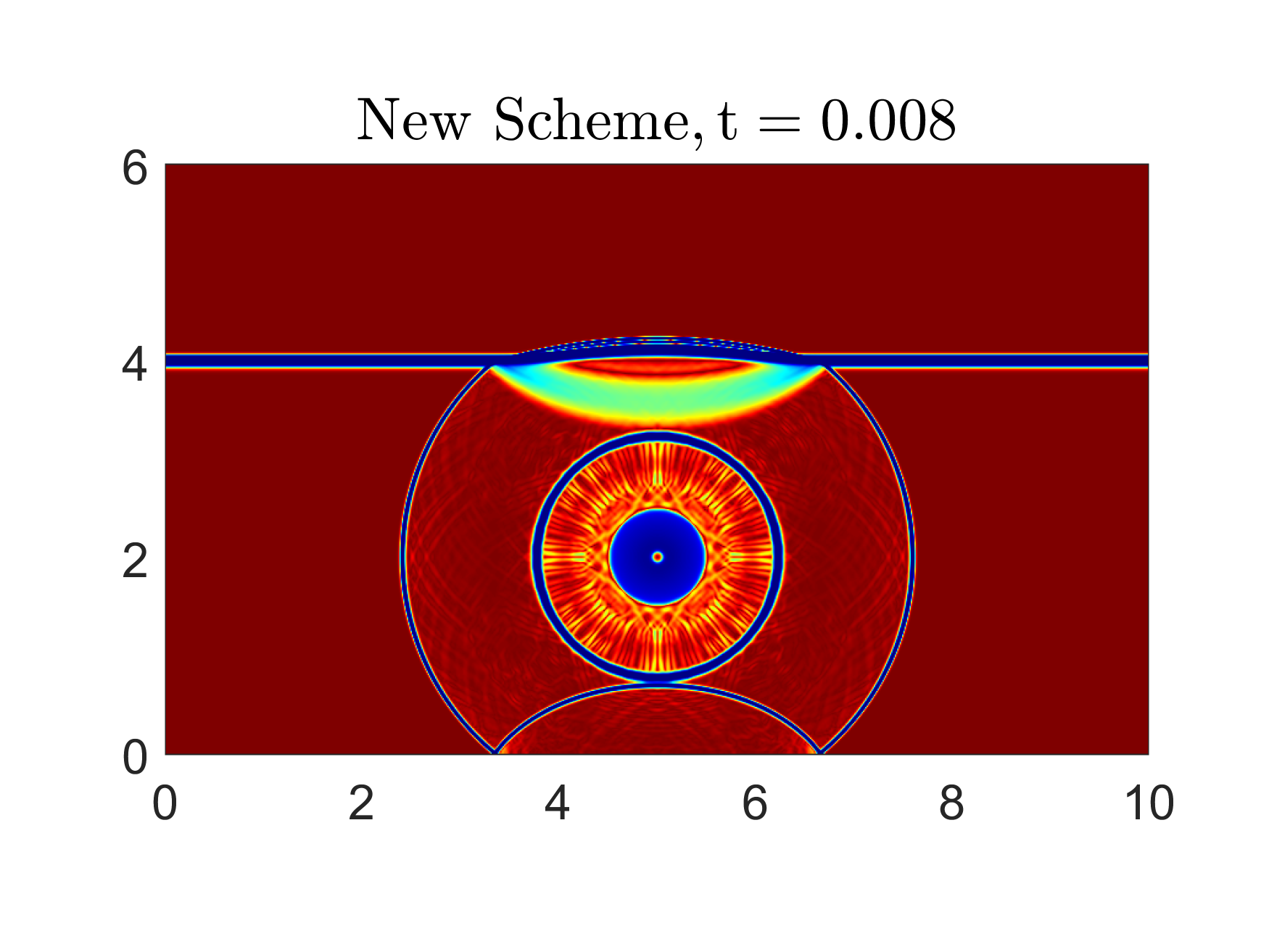}\hspace{1cm}
            \includegraphics[trim=1.2cm 1.5cm 1.1cm 1.1cm, clip, width=7cm]{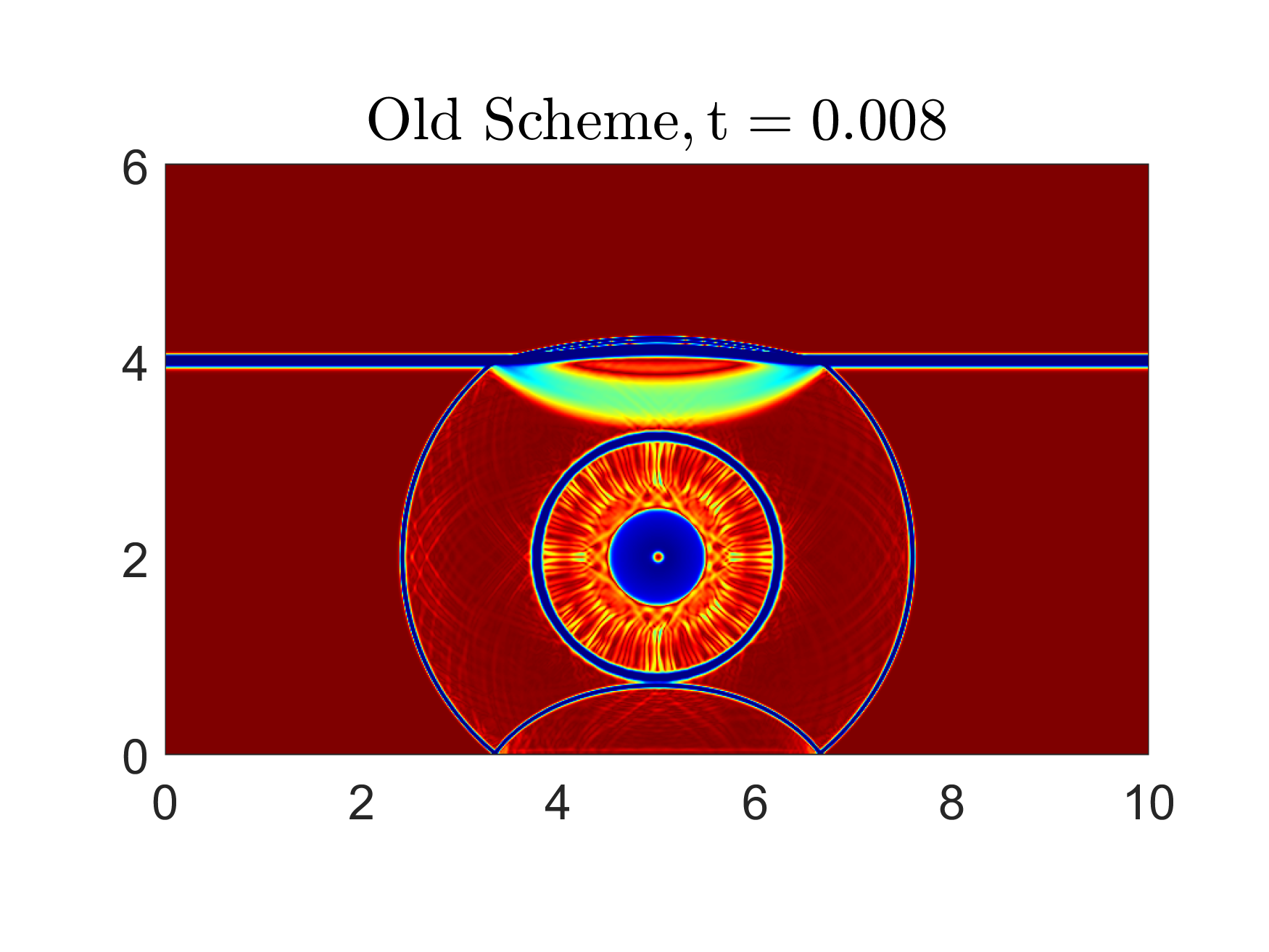}}
\vskip10pt
\centerline{\includegraphics[trim=1.2cm 1.5cm 1.1cm 1.1cm, clip, width=7cm]{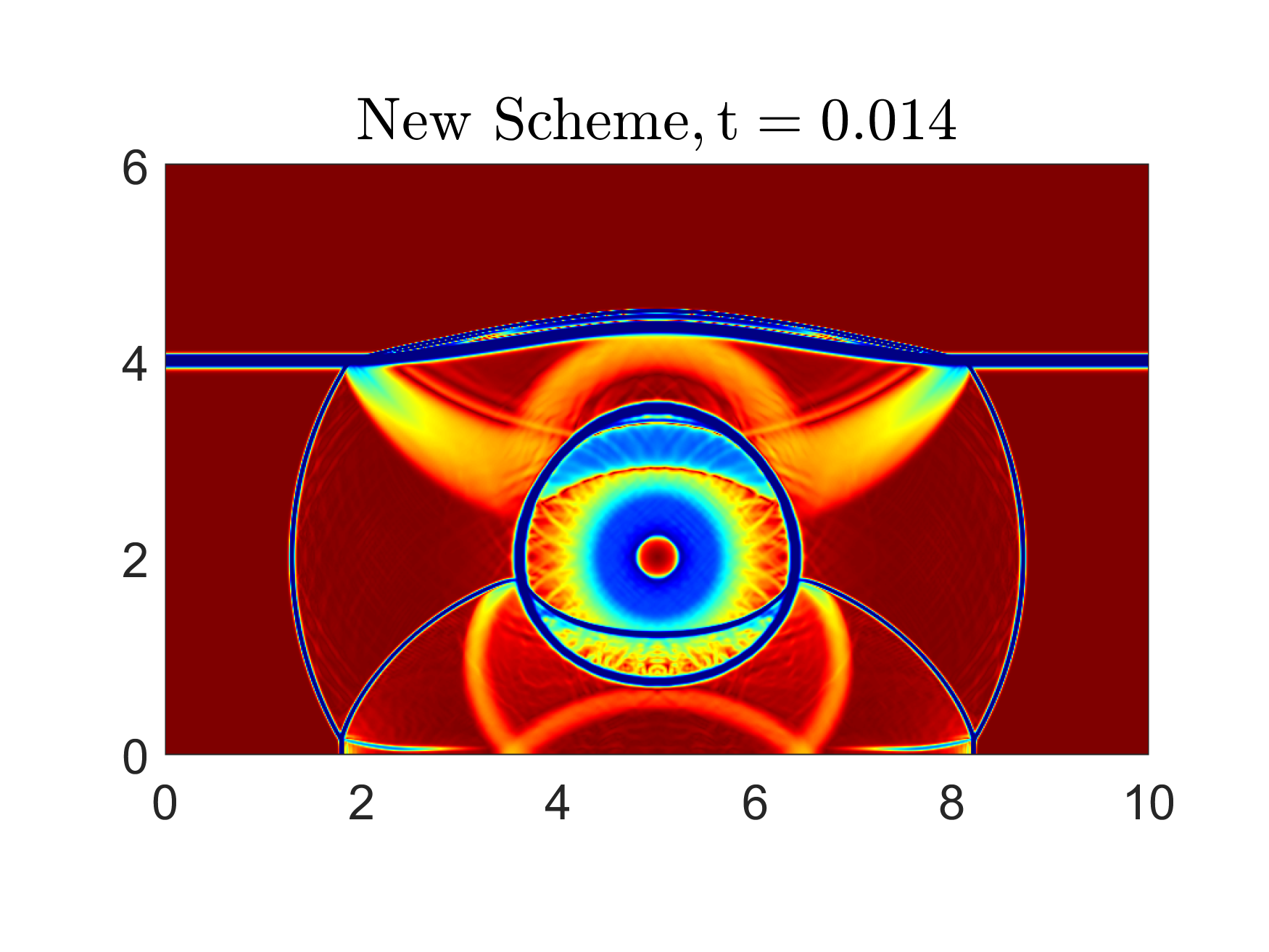}\hspace{1cm}
            \includegraphics[trim=1.2cm 1.5cm 1.1cm 1.1cm, clip, width=7cm]{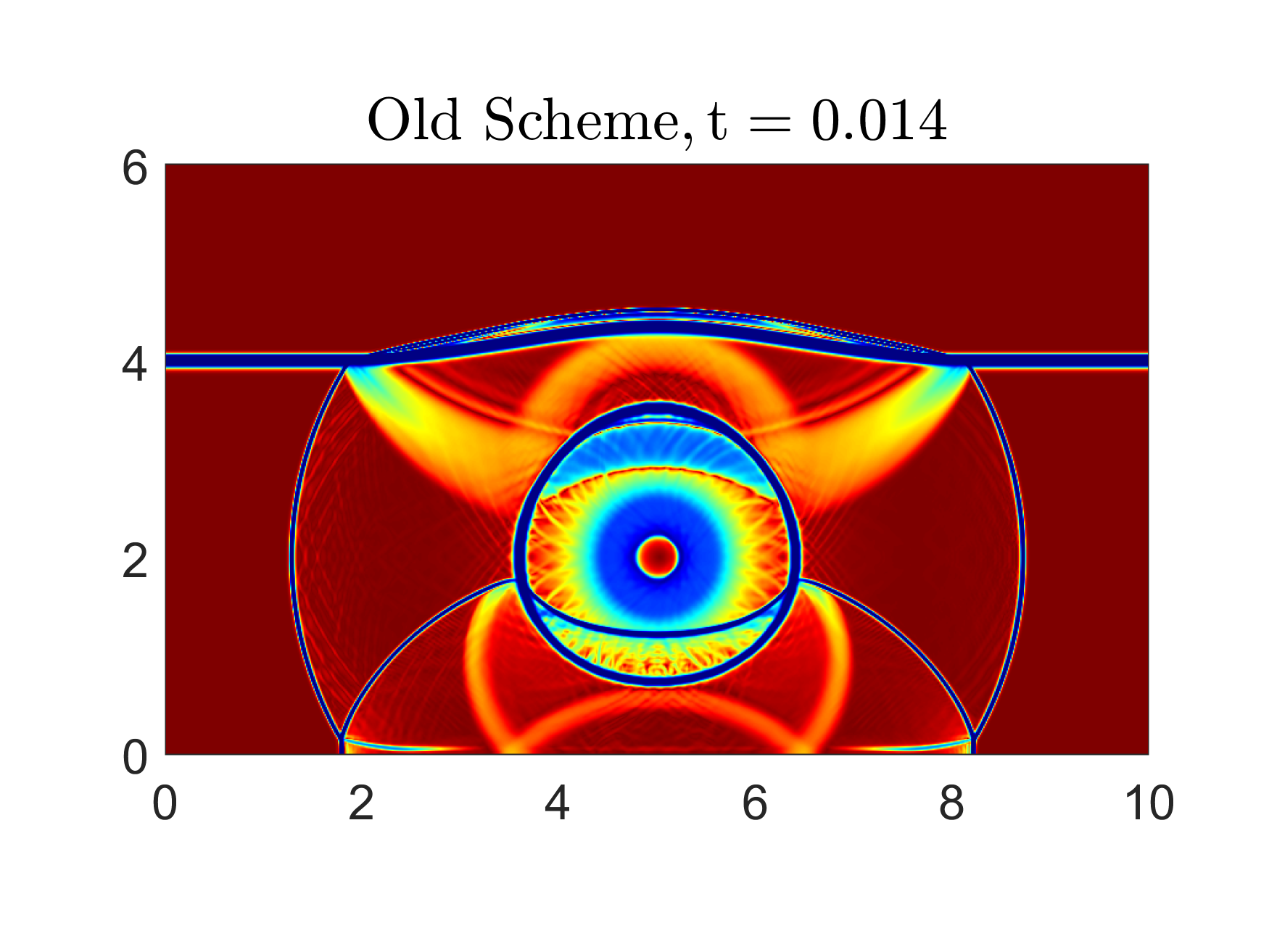}}
\vskip10pt
\centerline{\includegraphics[trim=1.2cm 1.5cm 1.1cm 1.1cm, clip, width=7cm]{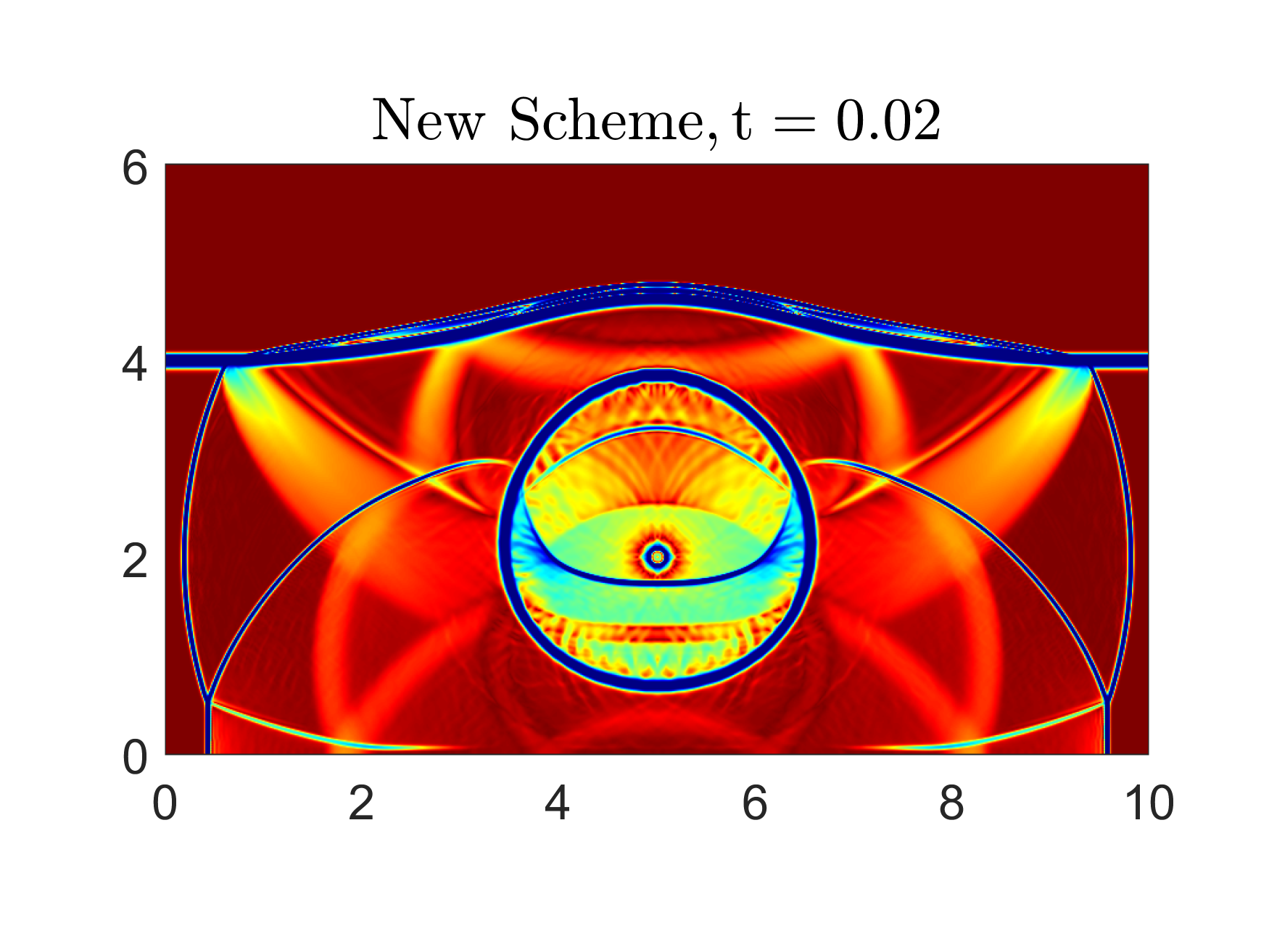}\hspace{1cm}
            \includegraphics[trim=1.2cm 1.5cm 1.1cm 1.1cm, clip, width=7cm]{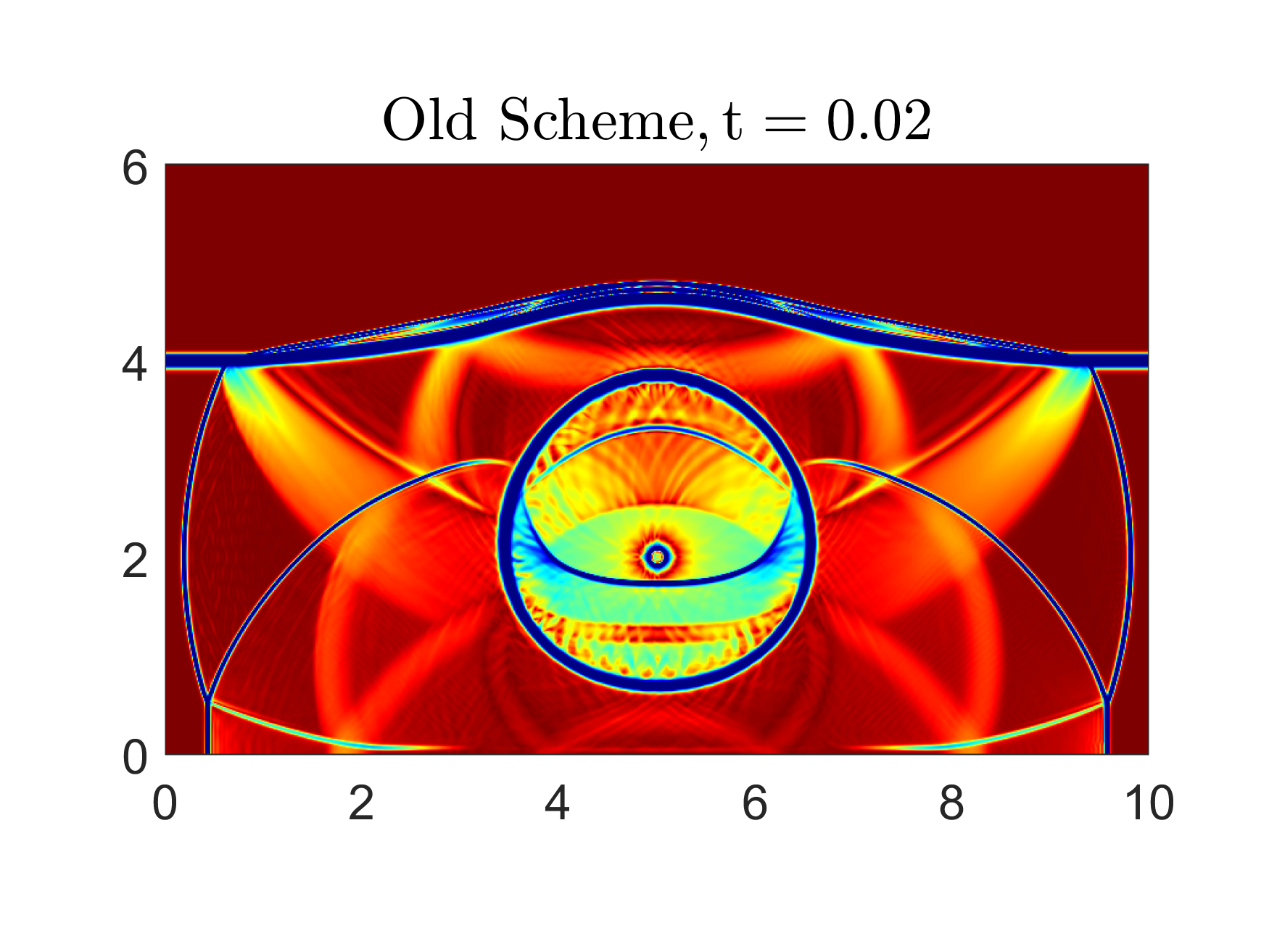}}
\caption{\sf Example 8: Density computed by the New (left) and Old (right) Schemes.\label{fig13}}
\end{figure}

We compute the numerical solutions until the final time $t=0.02$ on a uniform mesh with $\dx=\dy=1/80$ by the studied New and Old Schemes.
In Figure \ref{fig13}, we present time snapshots of the obtained results, where one can see that the results computed by the New and Old
Schemes are almost identical. 

We also measure the CPU times consumed by the New and Old Schemes. The obtained results are presented in Table \ref{tab_ex8}, where one can
see that the CPU time consumed by the Old Scheme is about $6.4\%$ larger than the CPU time consumed by the New Scheme. This finding
suggests that applying the new A-WENO schemes to nonconservative hyperbolic systems may yield greater efficiency improvements compared to
those observed in hyperbolic systems of conservation laws. Furthermore, it is important to note that, in contrast to the scalar equation
considered in Example 3, the efficiency gain is now smaller since due to the application of the 2-D local characteristic decomposition in
WENO-Z reconstructions (see, e.g., \cite{CKX24,CHK25}), which increased the computational cost for other components of the A-WENO scheme,
thereby diminishing the relative computational cost associated with the correction terms. 
\begin{table}[ht!]
\centering
\begin{tabular}{|c|c|c|}
\hline
Old Scheme&New Scheme&Ratio\\
\hline
4028.75&3785.95&1.064\\
\hline 
\end{tabular}
\caption{\sf Example 8: The CPU times (in seconds) consumed by the studied schemes.\label{tab_ex8}}
\end{table}

\section{Conclusions}
In this paper, we have developed new more efficient fifth-order A-WENO schemes for both conservative and nonconservative nonlinear
hyperbolic systems. A higher efficiency is achieved by utilizing the numerical fluxes in the computation of high-order correction terms. We
have performed a careful numerical study of the proposed schemes and demonstrated that they may be up to 15\% more efficient than the
existing A-WENO schemes and at the same time, neither accuracy of capturing smooth parts of the solutions nor the resolution of
discontinuous parts of the solution are affected by switching to the new approach.

\begin{DA}
\paragraph{Funding.} The work of S. Chu was supported in part by the DFG (German Research Foundation) through HE5386/19-3, 27-1. The work of
A. Kurganov was supported in part by NSFC grants 12171226 and W2431004.

\paragraph{Conflicts of interest.} On behalf of all authors, the corresponding author states that there is no conflict of interest.

\paragraph{Data and software availability.} The data that support the findings of this study and FORTRAN codes developed by the authors and
used to obtain all of the presented numerical results are available from the corresponding author upon reasonable request.
\end{DA}

\bibliography{refs}
\bibliographystyle{siamnodash}
\end{document}